\documentclass{amsart}
\usepackage{amsmath,amsthm,amssymb,amscd,stmaryrd}
\usepackage[mathscr]{eucal}

\numberwithin{equation}{section}
\numberwithin{figure}{section}

\theoremstyle{plain}
\newtheorem{theorem}[equation]{Theorem}

\newtheorem{lemma}[equation]{Lemma}
\newtheorem{proposition}[equation]{Proposition}

\theoremstyle{definition}
\newtheorem{definition}[equation]{Definition}

\theoremstyle{remark}
\newtheorem*{claim}{Claim}
\newtheorem*{remark}{Remark}

\newcommand{\Lra}{\Longrightarrow}
\newcommand{\Ra}{\Rightarrow}
\newcommand{\lra}{\longrightarrow}
\newcommand{\ra}{\rightarrow}
\newcommand{\up}[1]{\,\mathrm{#1}\,}
\newcommand{\be}{\begin{enumerate}}
\newcommand{\ee}{\end{enumerate}}

\newcommand{\beq}{\begin{equation}}
\newcommand{\eq}{\end{equation}}
\newcommand{\enq}{\end{equation}}
\newcommand{\beqs}{\begin{equation*}}
\newcommand{\eqs}{\end{equation*}}
\newcommand{\enqs}{\end{equation*}}
\newcommand{\recip}[1]{\dfrac{1}{#1}}

\newcommand{\inv}[1]{{#1}^{-1}}
\newcommand{\f}{\mathbb{F}_p}
\newcommand{\ff}{\mathbb{F}_2}
\newcommand{\Z}{\mathbb{Z}}
\newcommand{\z}{\mathbb{Z}/(p)}
\newcommand{\zn}[1]{\mathbb{Z}/(p^{#1})}
\DeclareMathOperator{\Tens}{Tens}
\DeclareMathOperator{\Ext}{Ext}
\DeclareMathOperator{\Ideal}{Ideal}
\DeclareMathOperator{\cl}{class}

\DeclareMathOperator{\Hom}{Hom}

\DeclareMathOperator{\Ind}{Ind}

\DeclareMathOperator{\image}{im}
\newcommand{\ol}[1]{\bar{#1}}
\newcommand{\bl}{\bigl(}
\newcommand{\br}{\bigr)}

\newcommand{\lp}{\bigl(}
\newcommand{\rp}{\bigr)}
\newcommand{\lb}{\bigl[}
\newcommand{\rb}{\bigr]}

\newcommand{\p}{\mathscr{P}}

\newcommand{\pti}[1]{\widetilde{\mathscr{P}}^{ {}^{\scriptstyle{#1}} }}
\newcommand{\bp}{\widetilde{\beta\mathscr{P}}^0}
\newcommand{\bpti}[1]{\widetilde{\beta\mathscr{P}}^{#1}} 
\newcommand{\bpit}[1]{\widetilde{\beta\mathscr{P}}^{#1}}

\newcommand{\sqit}[1]{\widetilde{Sq}^{#1}}
\newcommand{\sq}{\widetilde{Sq}^{1}}

\newcommand{\ot}{\otimes}
\newcommand{\sss}[1]{\scriptscriptstyle{#1}}
\newcommand{\sign}[1]{(-1)^{#1}}
\newcommand{\eo}{E^{\scriptscriptstyle 0}}

\DeclareMathOperator{\gr}{gr.}
\renewcommand{\L}{\mathscr{L}}
\newcommand{\V}{\mathscr{V}}
\newcommand{\VL}{\mathscr{VL}}
\newcommand{\vl}{\VL}

\begin{document}

\author{Justin Mauger}
\address{Whittier College\\
Whittier, CA 90608}
 \email{jmauger@whittier.edu}

\date{\today}

\title{On the Cohomology of the Lie Algebra Arising from the Lower Central Series of a $p$-Group.}

\begin{abstract}
We study the cohomology $H^*(A)=\Ext_A^*(k,k)$ of a locally finite, connected,
cocommutative Hopf algebra $A$ over $k=\f$. Specifically, we are interested in those
algebras $A$ for which $H^*(A)$ is generated as an algebra by $H^1(A)$ and $H^2(A)$.  We
shall call such algebras \emph{semi-Koszul}.  Given a central extension of Hopf algebras
$F\ra A\ra B$ with $F$ monogenic and $B$ semi-Koszul, we use the Cartan-Eilenberg
spectral sequence and algebraic Steenrod operations to determine conditions for $A$ to be
semi-Koszul.  Special attention is given to the case in which $A$ is the restricted
universal enveloping algebra of the Lie algebra obtained from the mod-$p$ lower central
series of a $p$-group.  We show that the algebras arising in this way from extensions by
$\z$ of an abelian $p$-group are semi-Koszul.  Explicit calculations are carried out for
algebras arising from rank 2 $p$-groups, and it is shown that these are all semi-Koszul
for $p\geq 5$.
\end{abstract}

\maketitle

\section{Introduction}

A \emph{quadratic} algebra over a field $k$ is an associative algebra of the form
$A=T(V)\big/(R)$ where $(R)$ is the two-sided ideal of relations generated by some
$R\subset V\ot V$, and $V$ is a locally finite graded vector space.
In \cite{Priddy}, Priddy defines a quadratic algebra $A$ to be a \emph{Koszul algebra} if
$H^*(A)=Ext_A^*(k,k)$ is generated as an algebra by $Ext_A^1(k,k)$. Basic examples
include tensor, polynomial, and exterior algebras.

When $A$ is Koszul, Priddy constructs projective resolutions over $A$ which are ``small''
subcomplexes of the bar resolution. He also shows that if $H^*(A)$ is itself Koszul, then
$A\simeq H^*(H^*(A))$ (Koszul duality). Equivalent definitions of Koszulity include:
\begin{enumerate}
\item A certain complex associated to $A$, the \emph{Koszul} complex, is acyclic
\cite[Theorem 1.2]{Lofwall}.
\item There exists a graded projective resolution
$$\dots\rightarrow P^2\rightarrow P^1\rightarrow P^0\rightarrow k\rightarrow 0$$
of $k$ such that $P^i$ is generated by its component of degree $i$, $P^i=A\cdot P^i_i$
\cite[Definition 1.2.1]{BGS}. Equivalently, $\Ext_A^{i,j}(k,k)=0$ if $i\neq j$.
\end{enumerate}

Now assume $A=T(V)/(R)$ is a more general algebra which is not necessarily quadratic,
i.e. $R\subset \sum_{i\geq 1} V^{\ot i}$. We are interested in the case when $H^*(A)$ is
generated as an algebra by $H^1(A)$ \emph{and} $H^2(A)$. Following Priddy, we shall call
such algebras \emph{semi-Koszul}. The quintessential example is the truncated polynomial
algebra $\f[x]/(x^{p})$, whose cohomology is generated by a one-dimensional exterior
class $z$ and a two dimensional polynomial class $u$. In this example, $u=-\bp (z)$. In
general, it will often be the case that the generators from $H^2$ are related to those of
$H^1$ by some algebraic Bockstein operation or Massey product. Since we allow $A$ to not
be quadratic, none of the results cited above for Koszul algebras necessarily hold for
semi-Koszul algebras.
It is unknown at this time whether there are any interesting relationships between $A$
and $H^*(H^*A)$, or anything like a semi-Koszul resolution.

We will focus on locally finite, connected, cocommutative Hopf algebras over $k=\f$.
Consequently, we will be able to use algebraic Steenrod operations. In Section \ref{main
theorem chapter}, we analyze extensions of Hopf algebras $F\ra A\ra B$, where
$F=\f[x]/(x^p)$ is central in $A$. There is a first quadrant spectral sequence converging
to $H^*(A)$ with $E_2^{**}\simeq H^*(B)\ot H^*(F)$.  The behavior of this spectral
sequence is determined by the differential $d_2(1\ot z)=\mu\in H^2(B)$, where $z\in
H^1(F)$. We prove in Theorem \ref{main theorem} that $A$ is semi-Koszul provided that $B$
is semi-Koszul, $\bp(\mu)=0$ and the annihilator ideal of $\mu$ is either zero or
generated in dimension one.

Section \ref{filtration} delves into Hopf algebras arising from $p$-groups. For a finite
$p$-group $G$, the mod-$p$ lower central series filtration $\{\Gamma^i G\}$ has as
associated graded object a restricted Lie algebra $\L G$. The universal restricted
enveloping algebra $\VL G$ is isomorphic to the associated graded algebra
$\eo (\f G)$  of the group ring $\f G$ with respect to the augmentation ideal filtration
\cite{Quillen}. Working with Lie algebras is helpful because many of the algebraic
Steenrod operations vanish there, in particular $\pti{0}$ \cite{Priddy-primary}. In
Section \ref{discussion}, we explore central extensions of $p$-groups $H\ra G\ra Q$,
where $H\simeq\z$.  Using Theorem \ref{main theorem}, we show that $\VL G$ is semi-Koszul
if $Q$ is abelian (Theorem \ref{base is abelian}). In Section \ref{vistas}, we prove a
corresponding statement for Hopf algebras which do not necessarily arise from $p$-groups.
Given a central extension of Hopf algebras $F\ra A\ra B$ with $F= \f[x]/(x^p)$, we show
that $\eo A$ is semi-Koszul provided that $B$ is commutative.

A rank two $p$-group is a finite $p$-group whose largest elementary abelian subgroup has
rank two. In Section \ref{rank 2} we show that $\VL G$ is semi-Koszul for all rank two
$p$-groups $G$ with $p\geq 5$.  One would hope to show that this holds for any finite
$p$-group, but such is not the case.  A counterexample is given in \ref{counterexample}.

Section \ref{algebraic preliminaries} introduces notation and covers the necessary
algebraic preliminaries such as algebraic Steenrod operations and the Cartan-Eilenberg
spectral sequence associated with a central extension of algebras. We relegate to the
Appendix definitions of the cobar complex $C^*(A)$, Massey products in $H^*(A)$ and
compuations in the cohomology of $\f[x]\big/(x^{p^n})$.


\section{Algebraic Preliminaries}\label{algebraic preliminaries}


In this section we introduce notation and conventions.  For more details on the
cohomology of algebras, the reader is referred to Adams \cite[Chapter 2]{Adams}.  We also
summarize algebraic Steenrod operations, the Cartan-Eilenberg spectral sequence, and Lie
algebras.

Let $k=\f$ be the field of $p$ elements. By an \emph{algebra} $A$ we mean a
non-negatively graded, associative $k$-algebra with product $\mu:A\ot A\ra A$, unit $\eta
: k\ra A$ and augmentation $\epsilon : A\ra k$. We assume $A$ is locally finite, that is
each $A_n$ is finite dimensional over $k$. We also assume $A$ is \emph{connected}, i.e.
$A_0=k$. Let $I=I(A)=\ker\epsilon$ denote the augmentation ideal. For $A$ bigraded
algebra

Tensor products are taken over $k$, unless otherwise noted. The tensor product of two
algebras $A$ and $B$ is defined by $(A\ot B)_n=$ $\sum_{i+j=n}A_i\ot B_j$, with
multiplication $(a\ot b)(x\ot y)=\sign{\deg b\, \deg x}ax\ot by$.  If the algebras in
question are bigraded, we define $(A\ot B)_{n,m}= \sum_{\genfrac{}{}{0pt}{}{i+s=n}{j+t=m}}   A_{i,j}\ot
B_{s,t}$, and the total degree is used in the multiplication formula. For $x\in A_{s,t}$,
$x$ has \emph{homological} degree (or dimension) $s$, \emph{internal} degree $t$ and
\emph{total} degree $s+t$. By \emph{commutativity} of $A$ we shall mean
\emph{graded-commutativity}, also called \emph{anticommutativity}: $ab=(-1)^{\deg a\,\deg
b}ba$.

A \emph{Hopf algebra} is an algebra with \emph{coproduct}, or \emph{diagonal,} $\psi:A\ra
A\ot A$ which is an algebra map.  We shall assume that $\psi$ is (graded) coassociative
and cocommutative. Since the algebra $A$ is assumed to be connected, it has an antipode
$\sigma:A\lra A$, but we shall not make use of it. See \cite{Milnor-Moore} for more
details on Hopf algebras and their structure.


Let $\Ext_A^i(-,k)$ denote the $i$-th right derived functor of $\Hom_A (-,k)$ in the
category of graded left $A$-modules.  The \emph{cohomology} $H^*(A)$ of $A$ is defined to
be $\Ext_A^*(k,k)$, where $k$ is given an $A$-module structure by the augmentation map
$\epsilon$. The Yoneda product makes this a graded algebra \cite{Adams}. 
We recall some basic facts about the cohomology of algebras.

\begin{enumerate}
\item Let $A,B$ be locally finite algebras. Then $H^*(A\ot B)\simeq H^*(A)\ot
H^*(B)$.
\item If $A$ is a connected cocommutative
Hopf algebra with unit, then $H^*(A)$ is anticommutative.
\end{enumerate}

We introduce the following definitions:
\begin{definition}\mbox{}

\begin{enumerate}
\item A graded algebra $C$ is called \emph{bigenerated} if it is generated as an
algebra by elements of homological degree 1 and 2.
\item An algebra $A$ is \emph{semi-Koszul} if $C=H^*(A)$ is bigenerated.
\end{enumerate}
\end{definition}

We shall make use of a reindexed form of \emph{algebraic Steenrod operations}. In
\cite[Theorem 11.8]{May}, May defines operations $\p^i :H^{s,t}(A)\ra
H^{s+(2i-t)(p-1),\,pt}(A)$ for $p\geq 3$.  These operations are zero if $2i<t$ or
$2i>s+t$. We would like to reindex these so that $\p^i$ raises homological degree by
$2i(p-1)$ and is trivial for $i<0$, $2i>s$. (May denotes these reindexed operations by
$\pti{i}$.) This regrading has the effect of making $\pti{0}$ the first possibly
non-trivial operation on $H^*(A)$. It also eliminates all operations on $H^{s,t}(A)$ for
$p>2$ when $t$ is odd. This won't be of concern to us as most of our algebras will be
concentrated in even degrees.

\begin{theorem}\label{steenrod}
 Let $A$ be a cocommutative Hopf algebra over $\f$.  There exist natural
homomorphisms
\begin{enumerate}
\item $\sqit{i}:H^{s,t}(A)\lra H^{s+i,2t}(A)$ for $p=2$,
\item  $\pti{i}:H^{s,2t}(A)\lra H^{s+2i(p-1),2pt}(A)$ and\\
       $\bpit{i}:H^{s,2t}(A)\lra H^{s+2i(p-1)+1,2pt}(A)$ for $p\geq 3$
\end{enumerate}
 with the following properties:
\begin{enumerate}
\item $\sqit{i}=0$ if $i<0$ or $i>s$,\\
      $\pti{i}=0$ if $i<0$ or $2i>s$,\\
      $\bpti{i}=0$ if $i<0$ or $2i\geq s$;
\item $\sqit{i}(x)=x^2$ if $i=s$,\\
      $\pti{i}(x)=x^p$ if $2i=s$.
\item \emph{Cartan formulae}\\
$\sqit{i}(xy)=\sum_{j=0}^i \sqit{j}(x)\sqit{i-j}(y)$,\\
$\pti{i}(xy)=\sum_{j=0}^i \pti{j}(x) \pti{i-j}(y)$,\\
  $\bpit{i}(xy)=\sum_{j=0}^i \bl\bpit{j}(x)\pti{i-j}(y)+\pti{j}(x)\bpit{i-j}(y)\br$;
\item \emph{Adem relations}
\begin{enumerate}
\item If $a<2b$,\\
 $\sqit{a}\sqit{b}=\sum_i \binom{b-i-1}{a-2i}\sqit{a+b-i}\, \sqit{i},$
\item If $a<pb$,\\
 $\widetilde{\beta^\epsilon \mathscr{P}}^{ \scriptstyle{a}} \pti{b}
=\sum_i(-1)^{a+i}\binom{(b-i)(p-1)-1}{a-pi} \widetilde{\beta^{\scriptstyle{\epsilon}}
\mathscr{P}}^{ {\scriptstyle{a+b-i}}}\, \pti{i}$,
\item If $a\leq pb$,\\
$\widetilde{\beta^\epsilon \mathscr{P}}^{ \scriptstyle{a}} \widetilde{\mathscr{P}}^{
{}^{\scriptstyle{b}} }=(1-\epsilon)\sum_i (-1)^{a+i}
\binom{(b-i)(p-1)-1}{a-pi} \widetilde{\beta\p}^{\scriptstyle{a+b-i}} \,\pti{i}$ \\
\phantom{${\beta^\epsilon \mathscr{P}}^{ \scriptstyle{a}} \widetilde{\mathscr{P}}^{
{}^{\scriptstyle{b}} }=$} $- \sum_i (-1)^{a+i}\binom{(b-i)(p-1)-1}{a-pi-1}
\widetilde{\beta^\epsilon\p }^{ {\scriptstyle{a+b-i}}}\, \bpit{i},$\\
\end{enumerate}
where $\epsilon=0$ or $1$, the binomial coefficients are taken modulo $p$ and, by abuse
of notation, $\widetilde{\beta^{\scriptstyle{0}}\p}^{ {\scriptstyle{s}} }=\pti{s}$ and
$\widetilde{\beta^{\scriptstyle{1}}\p}^{  {\scriptstyle{s}} }=\bpit{s}$.
\end{enumerate}
\end{theorem}

\begin{remark}
One can actually define operations $\pti{i/2}$, $\bpti{i/2}$ on elements of odd internal
degree such that the Cartan formulae and Adem relations still make sense.  See
\cite[Appendix 1.5]{Ravenel} for details.
\end{remark}


A \emph{decreasing filtration} of $A$ is a sequence of $k$-modules $A=F^0 A\supset F^1
A\supset F^2 A\supset\dotsb$ such that $F^s A\cdot F^t A\subset F^{s+t}A$. Let $\rho_s$
denote the natural map $F^s A\twoheadrightarrow F^{s/s+1}A=F^sA\big/F^{s+1}A$. We will
sometimes suppress the $A$ when it is clear from context. We will say that elements in
$F^{s+i}$ are of \emph{lower} filtration level than elements in $F^s$, even though the
indices are greater. Let the \emph{filtration degree} of an element $x\in A$ be the
largest $i$ such that $x\in F^i A$ (if such an $i$ exists).

Let $F\subset A$ be a subalgebra. We say $F$ is \emph{normal} in $A$ if $I(F)\cdot
A=A\cdot I(F)$. In this case, there is a well defined  quotient $B= A\big/ A\cdot I(F)$
denoted by $A//F$.  In the case that $F$ is a sub-\emph{Hopf} algebra, $B$ is also a Hopf
algebra. We say $F$ is \emph{central} in $A$ if each element of $F$ anticommutes with
each element of $A$. By an \emph{extension} or \emph{exact sequence} $F\ra A\ra B$ of
algebras we shall mean that $F$ is normal in $A$ and $B\simeq A//F$.
\begin{theorem}\label{Cartan-Eilenberg} Let $F \overset{i}{\ra} A \overset{\pi}{\ra} B$ be a
central extension of graded, connected algebras.  There is a first quadrant spectral sequence converging to $H^*(A)$ with
\be
\item $E_2^{s,t}\simeq H^s(B)\otimes H^t(F)$.
\item $E_\infty^{s,t} \simeq F^sH^{s+t}(A)\big/ F^{s+1}H^{s+t}(A)$ where $\{F^s\}$ is a
decreasing filtration of $H^*(A)$.
\item The ring structure of $H^s(B)\ot H^t(F)$ is defined by \\
$(x\ot y)(z\ot w)=\sign{(t+t')(s+s')} xz\ot yw$ where $y\in H^{t,t'}(F)$, $z\in
H^{s,s'}(B)$.
\item $d_r:E_r^{s,t}\ra E_r^{s+r,t-r+1}$ is a graded derivation, meaning \\ 
$d_r(ab)=d_r(a)\cdot b +\sign{\deg a}a \cdot d_r(b)$.
\item The map $H^s(B)\simeq E_2^{s,0}\twoheadrightarrow E_\infty^{s,0}\simeq F^sH^s(A)
\rightarrowtail H^s(A)$ corresponds to the induced map $\pi^*:H^*(B)\ra H^*(A)$. \\
      The map $H^t(A)=F^0H^t(A)\simeq E_\infty^{0,t}\rightarrowtail E_2^{0,t}\simeq H^t(F)$
corresponds to the induced map $i^*:H^*(A)\ra H^*(F)$.
\item \emph{Kudo transgression theorem.} Assume, in addition, that
$F \overset{i}{\ra} A \overset{\pi}{\ra} B$
is an extension of Hopf algebras over $\f$.
 The \emph{transgression} $d_r:E_r^{0,r-1}\ra E_r^{r,0}$ anticommutes with Steenrod operations:
if $x\in E_r^{0,r-1}$ and $\theta$ is a Steenrod operation of homological
 degree $i$, then $\theta(x)$ survives until $E_{r+i}^{0,r+i-1}$ and
 $d_{r+i}\bl\theta (x)\br=\sign{i}\theta\bl d_r(x)\br$.
\ee
\end{theorem}

Since the algebras $F$, $A$ and $B$ are graded, the spectral sequence is actually
trigraded, with the third grading being the internal degree.  We have suppressed this
third grading, as it is preserved by the differentials $d_r$. This spectral sequence was
used by Adams \cite[2.3.1]{Adams}. Ravenel \cite[A1.3.14]{Ravenel} calls it the
\emph{Cartan-Eilenberg} spectral sequence.  For (6) see \cite[Proposition
2.3]{Wilkerson}.


The Lie algebras we shall work with are all graded over either $\Z$ or $\f$.  A
\emph{restricted} Lie algebra is a Lie algebra over $\f$ with restriction $\xi$.

For $L$ a restricted Lie algebra, let $\Tens(L)$ be the tensor algebra on $L$, and let
$J\subset \Tens(L)$ be the two sided ideal generated by elements of the form $x\ot
y-(-1)^{\deg x\,\deg y}y\ot x-[x,y]$ and $x^{\ot p}-\xi(x)$ for $x,y\in L$.
\begin{definition}\label{VL}
The \emph{universal restricted enveloping algebra} of $L$ is defined to be $\V
L=\Tens(L)/J$.
\end{definition}
We can put a Hopf algebra structure on $\V L$ by setting $\psi(x)=x\ot 1+1\ot x$,
$\epsilon(x)=0$ for $x\in L$. The cohomology of a Lie algebra $L$ is defined to be
$H^*(\V L)$.
\begin{theorem}[{ \cite[5.3]{Priddy-primary}, \cite[Theorem 8.5]{May} }]
Let $L$ be a Lie algebra over $\f$. Then the Steenrod operations $\pti{0}$ and $\sqit{0}$
are identically zero on $H^*(\V L)$.
\end{theorem}


\section{An Extension Theorem}\label{main theorem chapter}

Let $A$ be a cocommutative  Hopf algebra over $\f$. In this section, we will find
sufficient conditions for $H^*(A)$ to be bigenerated. Let $z\in A$ be a central primitive
element of height $p$. This means $z$ anticommutes with every element of $A$, $z^p=0$ and
$z^i\neq 0$ if $i<p$. Let $F=\f[z]\big/(z^p)$ be the central sub-Hopf algebra of $A$
generated by $z$. Let $B=A//F$.  We have the central extension of  Hopf algebras
\beq\label{main theorem Hopf algebra extension}
F\lra A\lra B
\enq
with corresponding spectral sequence
\begin{equation}\label{main theorem Hopf algebra spectral sequence}
E_2^{s,t}\simeq H^s(B)\ot H^t(F)\Rightarrow H^{s+t}(A).
\end{equation}
Let $n$ be the internal degree of $z\in F$.  Note that when $p$ is odd, $n$ must be even.
We recall that the cohomology of the fiber $F$ is
\begin{align*}
H^*(F) &\simeq \Lambda[\zeta]\ot \f[\epsilon],& \zeta&\in H^{1,n},
\,\epsilon\in H^{2,np},& \text{with }\bp(\zeta)&= -\epsilon & p \text{ odd,}\\
H^*(F) &\simeq \ff[\zeta],& \zeta&\in H^{1,n},& \text{with }\sqit{1}(\zeta)&=\zeta^2 &
p=2.
\end{align*}
See Theorem \ref{non-zero Steenrod operations} for more details on the Steenrod
operations.

\begin{theorem}[Main Theorem] \label{main theorem} Let $F,A,B$ be as in
\eqref{main theorem Hopf algebra extension}. In the spectral sequence \eqref{main theorem
Hopf algebra spectral sequence}, let $d_2(\zeta)=\mu\in H^{2,n}(B)$. Assume
\be
\item $\bp(\mu)=0$ or $\sq{\mu}=0$.
\item The annihilator ideal of $\mu$ in $H^*(B)$ is either zero or is generated in dimension
one.
\item $B$ is semi-Koszul.
\ee
Then $A$ is semi-Koszul.
\end{theorem}

\noindent\textit{Proof.} By the transgression theorem \ref{Cartan-Eilenberg}(6),
\begin{align*}
d_3(\epsilon)&=-\bp(\mu), \quad p \text{ odd},\\
d_3(\zeta^2) &=\sq (\mu), \quad\, p=2.
\end{align*}
If $\bp(\mu)=0$ or $\sq(\mu)=0$ then
\begin{align*}
&d_{2 p+1}(\epsilon^p) =\pti{1}\bl d_3(\epsilon)\br=0 \quad\text{ or}\\
&d_5(\zeta^4) =\sqit{2}\bl d_3(\zeta^2)\br=0.
\end{align*}
Likewise for all other powers of $\epsilon$ or $\zeta$.  Thus in this case, the spectral
sequence \eqref{main theorem Hopf algebra spectral sequence} collapses at $E_3$.

In order for $H^*(A)$ to be bigenerated, we need $E_\infty$ to be bigenerated. Assume
$\bp(\mu)=0$ or $\sq(\mu)=0$. If $\mu$ is a zero divisor in $H^*(B)$, we may run into the
following problem. Let $\eta\in H^i(B)$ be in the annihilator ideal of $\mu$. Then
$$d_2(\eta\ot\zeta)=d_2(\eta)\ot\zeta +\sign{i}\eta \ot d_2(\zeta)=\pm\eta\ot\mu=0.$$
Nothing can hit $\eta\ot\zeta\in E_r^{i,1}$ because $d_r\equiv 0$ for $r\geq 3$, and
$d_2\equiv 0$ from the $E_2^{*,2}$ line. Thus $\eta\ot\zeta$ is a non-bounding permanent
cycle of dimension $i+1$. If $i=1$, then $\eta\ot\zeta\in E_\infty^{1,1}$ is a
two-dimensional generator of $E_\infty$.  If $i>1$ on the other hand, $\eta\ot\zeta$ is a
possibly indecomposable element in $E_\infty^{i,1}$ of homological dimension larger than
$2$. (It will be indecomposable if, for example, $\eta$ is indecomposable in $H^*(B)$.)
Consequently $E_\infty$---and thus $H^*(A)$---is not bigenerated.

Assume then that the annihilator ideal of $\mu$ in $H^*(B)$ is either zero or is
generated multiplicatively by elements of dimension one.
Assume also that $H^*(B)$ is bigenerated. Then $E_\infty$ has the following possible set
$S$ of multiplicative generators:
\be
\item $1\ot \epsilon$ for $p$ odd, $1\ot\zeta^2$ for $p=2$;
\item elements of the form $\eta\ot\zeta$, where $\eta\in H^1(B)$ and $\eta\mu=0$;
\item elements $\sigma\ot 1$, $\nu\ot 1$, where $\sigma\in H^1(B)$ and $\nu\in H^2(B)$;
\item $1\ot \zeta$ (if $d_2\equiv 0$).
\ee

All these elements have homological degree less than or equal to two, and thus $E_\infty$
 is  bigenerated.  Recall that
$H^*(A)$ has a decreasing filtration $F^*H^*(A)$, and that we denote the natural
projection $F^sH^n(A)\ra F^{s/s+1}H^n(A)\simeq E_\infty^{s,n-s}$ by $\rho_s$. Choose a
set of representatives $T\subset H^*(A)$ for each of the types of generators in $S$:
\be
\item $e\in F^0H^2(A)$;
\item $x\in F^1 H^2(A)$;
\item $\alpha\in F^1H^1(A)$, $v\in F^2 H^2(A)$;
\item $z\in F^0H^1(A)$.
\ee

We claim that $H^*(A)$ is generated as an algebra by $T$, from which it follows that $A$
is semi-Koszul. Let $S'$=$S\cup \{1\ot 1\}$, $T'=T\cup\{1\}$. The proof will proceed by
reverse induction on filtration degree $s$ in $H^n(A)$. Let $y\in F^nH^n(A)$. Then $y\in
H^n(B)$, which is bigenerated. Now let $y\in F^sH^n(A)$, where $s<n$, and assume
$F^iH^n(A)$ is bigenerated for $i>s$. Projecting $y$ to $F^{s/s+1}H^n(A)\simeq
E_\infty^{s,n-s}$, $\rho_s(y)$ can be written as a sum of products
$f_i\lambda_{i_1}\lambda_{i_2}\dotsm\lambda_{i_n}$ where $f_i\in k$ and $\lambda_{i_j}\in
S'$. So $y$ is a sum of products $f_i\tilde{\lambda}_{i_1}\dotsm \tilde{\lambda}_{i_n}$
and a term $\gamma$, where $\tilde{\lambda}_{i_j}\in T'$ are lifts of $\lambda_{i_j}$
under $\rho$, and $\gamma\in F^{s+1}H^n(A)$.  By induction, $\gamma$ is also a linear
combination of products of elements from $T'$. Thus $y$ is of the desired form. \qed


\section{Hopf algebras arising from $p$-groups}\label{filtration}

The mod-$p$ lower central series of a group $G$ gives rise to a restricted Lie algebra
$L$, whose universal restricted enveloping algebra $VL$ is a Hopf algebra.  Starting with
a central extension of $p$-groups $H\ra G\ra Q$, we take enveloping algebras of the
associated Lie algebras and get a central extension of Hopf algebras $F\ra A\ra B$. In
Section \ref{discussion} we show that $A$ is semi-Koszul if $Q$ is abelian and
$H\simeq\z$. In Section \ref{vistas}, we generalize to extensions of Hopf algebras $F\ra
A\ra B$ which do not necessarily arise from extensions of $p$-groups.  We show $\eo A$ is
semi-Koszul if $B$ is commutative.

\subsection{Definitions} Let $G$ be a group, and let $G=F^1G\supseteq
F^2G\supseteq \dotsb$ be a descending sequence of subgroups (also called a filtration)
with $F^iG\unrhd F^{i+1}G$, \,$F^iG\big/F^{i+1}G$ abelian. Let
$\mathrm{gr}.^FG=\bigoplus_i F^iG\big/F^{i+1}G$ be the \emph{associated graded} object,
and let ${\mathrm{gr}_i}^FG=F^iG\big/F^{i+1}G$.  When the filtration is understood, we
shall simply write $\gr G$ and $\mathrm{gr}_i G$. Recall the natural surjection $\rho_i
:F^iG\twoheadrightarrow \mathrm{gr}_i G$.

For $H,\, K$ subgroups of $G$, let $(H,\, K)$ be the subgroup generated by all
commutators of the form $(h,\, k)=\inv{h}\inv{k}hk$ for $h\in H,\, k\in K$. Let $H^n$ be
the subgroup generated by $h^n$ for $h\in H$. An \emph{$N$-sequence} is a normal series
$G=G_1\unrhd G_2\unrhd G_3\unrhd \dotsb$ such that $(G_i,\, G_j)\subset G_{i+j}$
\cite{Lazard}. The successive quotients $\mathrm{gr}_i G=G_i\big/ G_{i+1}$ of an
$N$-sequence are abelian, since $(G_i,G_i)\subset G_{2i}\subset G_{i+1}$. The associated
graded abelian group $\gr G=\bigoplus_n G_n\big/ G_{n+1}$ has the structure of a Lie
algebra over $\Z$, with bracket $\lb\rho_r x,\rho_s y\rb =\rho_{r+s}\bl(x,y)\br$. A
\emph{$p$-restricted} $N$-sequence is an $N$-sequence $\{G_n\}$ such that ${G_n}^p
\subset G_{pn}$.  The associated graded group $\gr G$ is then a vector space over $\f$.
We can define a restriction $\xi\bl\rho_r(x)\br=\rho_{pr}\bl x^p\br$ on $\gr G$ which
makes it into a restricted Lie algebra over $\f$.

The fastest descending $p$-restricted $N$-sequence of $G$ is the \emph{mod-$p$ lower
central series} $\bigl\{ \Gamma^iG \bigr\}$: (see \cite{Lazard}, \cite{Blackburn2})
\begin{equation}\label{filtration def}
 \Gamma^{2n-1}G=\Gamma^{2n}G =
 \langle (\gamma_1 , \dots , \gamma_r )^{p^u} \text{ such that } r\cdot p^u
 \geq n \rangle,
\end{equation}
where $(\gamma_1 , \dots , \gamma_r)=
 \bl\gamma_1,\bl \dotso \bl\gamma_{r-2},(\gamma_{r-1},\gamma_r)\br\dotso\br$.
The associated object $\mathrm{gr}.^\Gamma G$ is thus a restricted Lie algebra over $\f$.
Let $\L(G)=\L G=\mathrm{gr}.^{\Gamma}G$. It can be easily checked that $\mathscr{L}$ is a
functor from the category of groups to the category of restricted Lie algebras over $\f$.

The definition we have given here is slightly different from that of \cite{Quillen} or
\cite{Blackburn2}.  Our $\Gamma^{2n}$ equals their $\Gamma^n$. We have doubled the
filtration degree in order not to have to worry about signs for $p$ odd.

Let $kG$ be the group ring of $G$ over $k$, with augmentation $\epsilon(g)=1$ for $g\in
G$. A diagonal $\psi:kG\ra kG\ot kG$ is induced by the map $g\ra g\ot g$, giving $kG$ the
structure of a cocommutative Hopf algebra over $k$. Let $I=\ker\epsilon$ be the
augmentation ideal of $kG$.  Filter $kG$ as follows: $F^{2n-1}=F^{2n}=I^n$ for $n>0$. The
associated graded algebra $\eo (kG)=\mathrm{gr}.^FkG$ has the structure of a restricted
Lie algebra over $k$ inherited from the multiplication in $kG$. Specifically,
$[\rho_{r}x,\rho_{s}y]= \rho_{r+s}(xy-yx),\, \xi(\rho_{r}x)=\rho_{rp}(x^p)$. Since $\psi$
respects the filtration on $kG$, it induces a cocommutative diagonal on $\eo (kG)$.  It
is clear that $\eo (kG)$ is generated as an algebra by $\eo_2 (kG)=I/I^2$. These elements
are all primitive by reason of dimension.  Thus $\eo (kG)$ is primitively generated, and
by a theorem of Milnor-Moore \cite[Theorem 6.11]{Milnor-Moore}, $\eo (kG)$ is isomorphic
to the restricted universal enveloping algebra of its restricted Lie subalgebra
$L=\eo_2 kG$.

The map $\psi:\L G\ra\eo (kG)$ given by $\rho_r(x) \mapsto\rho_r(x-1)$ is a restricted
Lie algebra homomorphism. In fact, the image of $\psi$ lies in the primitive elements
$\mathcal{P}\bigl(\eo (kG)\bigr)$ of $\eo (kG)$.
\begin{theorem}[{\cite[Theorem 1.4]{Quillen}}]\label{Quillen's theorem} The induced map
 $$\L G\longrightarrow \mathcal{P}\bl\eo (kG)\br$$
is an isomorphism of restricted Lie algebras.
\end{theorem}
  Therefore,
\begin{equation}\label{previous equation}
\V(\L G) \simeq \V\bigl(\mathcal{P}(\eo (kG))\bigr) \simeq \eo (kG)
\end{equation}
as Hopf algebras over $k$.

Let $\VL G$ denote $\V(\mathrm{gr}.^{\Gamma}G)$. Studying the cohomology $H^*\bl\eo
(kG)\br$ over $k$ is thus equivalent to studying $H^*(\VL G)$. The advantage of working
with cohomology of Lie algebras is that the Steenrod operations $\sqit{0}\equiv 0$,
$\pti{0}\equiv 0$  there (\cite{Priddy-primary}, \cite[Theorem 8.5]{May}), which greatly
simplifies the situation.

We shall need the following result.
\begin{lemma}\label{VL for G abelian}
Let $G$ be a finite abelian $p$-group, then $\VL G$ is a semi-Koszul algebra.
\end{lemma}
\begin{proof}
 Let $G=\zn{n}$. It is clear that
\begin{equation*}
\VL G\simeq \eo (\f G)\simeq \f[x] \big/ \bl x^{p^n}\br
\end{equation*}
For a general finite $p$-group $G=\zn{n_1}\times\dots\times\zn{n_k}$, then
\begin{equation*}
\VL G\simeq \f[x_1,\dots, x_k]\big/\bl x_i^{p^{n_i}} \br
\end{equation*}
Thus
\begin{equation*}
H^*(\VL G)\simeq \Lambda[z_1,\dots ,z_k]\ot\f[e_1,\dots,e_k],
\end{equation*}
where $z_i\in H^{1,2}$,\, $e_i\in H^{2,2p^{n_i}}$.  If $p=2$ and some $n_i=1$, replace
$\Lambda[z_i]\ot\f[e_i]$ by $\ff[z_i]$.
\end{proof}

\subsection{An Extension Theorem}\label{discussion}

Let $G$ be a non-trivial finite $p$-group, and let $k=\f$.  There is an element $x$ of
order $p$ in the center of $G$. Let $H\simeq\z$ be the subgroup generated by $x$. We have
the central extension
\begin{equation}\label{group extension}
H\overset{f}{\lra} G\overset{g} {\lra}Q=G/H
\end{equation}
Our plan is to apply the $\mathscr{L}$ functor, and then the $\mathscr{V}$ functor, to
(\ref{group extension}) in the hope of obtaining an exact sequence of Hopf algebras of
the form \eqref{main theorem Hopf algebra extension}.

First, note that $g(\Gamma^n G)=\Gamma^n Q$.  Refilter $H$ by
$F^nH=\inv{f}\bl\Gamma^nG\br$.  According to Theorem 2.4 of \cite{Lazard},
the following is an exact sequence of graded Lie algebras over $\Z$.
\begin{equation}\label{seq}
\gr^FH\overset{\ol{f}}{\lra} \gr G \overset{\ol{g}}{\lra} \gr Q
\end{equation}
Now $\mathrm{gr}.^FH$ is generated by $\ol{x}\in F^{2n/2n+1}H$, where $2n$ is the
filtration degree of $f(x)$ in $G$. Since $\L H$ is generated by $\ol{x}\in
\Gamma^{2/3}H$, $\mathrm{gr}.^FH\simeq \L H$ as ungraded Lie algebras.  The sequence
\eqref{seq} is actually a sequence of graded restricted Lie algebras and---if we regrade
$\L H$ appropriately---is isomorphic to
\begin{equation}\label{Lie algebra extension}
\L H\lra \L G \lra \L Q.
\end{equation}

Now let
\begin{equation}
L_1\overset{i}{\lra} L_2 \overset{\pi}{\lra} L_3
\end{equation}
be an exact sequence of graded restricted Lie algebras over $\f$. Recall that $\V L_i$ is
defined to be the quotient of $\Tens(L_i)$ by some specific ideal $J_i$. It is clear that
$\pi(J_2)=J_3$, and that $\inv{i}(J_2)=J_1$. Therefore,
\begin{equation}
\Tens(L_1)\big/J_1 \lra \Tens(L_2)\big/J_2 \lra \Tens(L_3)\big/J_3
\end{equation}
is an exact sequence of Hopf algebras over $k$.

We have just proven
\begin{proposition} Let $H, G,$ and $Q$ be as in (\ref{group extension}).
Regrade $\VL H$ as in \eqref{Lie algebra extension}.
Then the following sequence of Hopf algebras is exact:
\beq\label{VL extension}
\VL H\lra \VL G\lra \VL Q
\enq
\end{proposition}

Since $\vl H\simeq \f[x]\big/(x^p)$, the central extension \eqref{VL extension} is of the
form $F\ra A\ra B$, as in \eqref{main theorem Hopf algebra extension}. Associated to this
extension is the Cartan-Eilenberg spectral sequence
 \eqref{main theorem Hopf algebra spectral sequence}:
\begin{equation}\label{algebra spectral sequence}
\begin{split}
E_2&\simeq H^*(\VL Q)\ot \Lambda[z]\ot\f[e]\Rightarrow H^*(\VL G), \quad p \text{
odd},\\
E_2&\simeq H^*(\VL Q)\ot \ff[z]\Rightarrow H^*(\VL G), \qquad\qquad p=2.
\end{split}
\end{equation}
Let us now examine conditions on the group $G$ which will assure that the
assumptions of Theorem \ref{main theorem} are satisfied.
\begin{theorem} In the spectral sequence \eqref{algebra spectral sequence}, assume
$\mu=d_2(z)$ is a sum of products of 1 dimensional elements.  That is, assume $\mu$ is
\emph{decomposable}. Assume also that Conditions 2 and 3 of Theorem \ref{main theorem}
hold. Then $\VL G$ is semi-Koszul.
\end{theorem}
\begin{proof} This result follows from the fact that
\begin{equation}\label{proof that bp}
\begin{aligned}\
\bp(\sigma\tau)&=\bp(\sigma)\pti{0}(\tau)+\pti{0}(\sigma)\bp(\tau)=0 \text{ and}\\
\sqit{1}(\sigma\tau)&=\sqit{1}(\sigma)\sqit{0}(\tau)+\sqit{0}(\sigma)\sqit{1}(\tau)=0,
\end{aligned}
\end{equation}
since $\pti{0}\equiv 0$ and $\sqit{0}\equiv 0$ on the cohomology of a Lie algebra.
Condition $1$ of Theorem \ref{main theorem} is thus satisfied.
\end{proof}
\begin{theorem}\label{base is abelian}
If $Q$ is abelian, then $\VL G$ is semi-Koszul.
\end{theorem}
\begin{proof}
By Corollary \ref{VL for G abelian}, $H^*(\VL Q)$ is isomorphic to a tensor product of
algebras of the form $\Lambda[z_i]\ot\f[e_i]$ (and $\ff[z]$ for $p=2$), so $H^2(\VL Q)$
is spanned by $\{z_iz_j, e_i\}$. Now $\bp(z_iz_j)=0=\sqit{1}(z_iz_j)$ by \eqref{proof
that bp} and $\bp(e_i)=0$ by reason of internal degree.  Thus Condition $1$ of Theorem
\ref{main theorem} is satisfied. Condition 2 is satisfied because the annihilator ideal
of $H^2(\VL Q)$ is generated by one-dimensional elements $\{z_i\}$. Finally, Condition 3
is satisfied since $H^*(\VL Q)$ is bigenerated.
\end{proof}

\subsection{Generalization to Hopf Algebras}\label{vistas}

The results in \S \ref{discussion} can be generalized somewhat to Hopf algebras.  Let $A$
be a Hopf algebra over $\f$.  Assume $A$ contains a central, primitive element $x$ of
height $p$. This means $x$ anticommutes with every element of $A$, $x^p=0$, and $x^i\neq
0$ if $i<p$.  If $p$ is odd, $x$ is forced to have even internal degree. Let
$M=\f[x]/(x^p)$ be the central sub-Hopf algebra of $A$ generated by $x$.

We have the central extension of graded Hopf algebras
\begin{equation}\label{general extension}
M \overset{f}{\lra} A \overset{g}{\lra} B.
\end{equation}
For an augmented algebra $C$, let $\{F^{2n-1}C=F^{2n}C=(IC)^n\}$ be the modified
augmentation ideal filtration, and let $\eo C=\mathrm{gr.}^FC$. Since
$\epsilon^B(g(y))=\epsilon^A(y)$, we see that $g(IA)\subset IB$, and more generally,
$g(IA^{\,n})\subset IB^{\,n}$ .  Thus $g$ respects the augmentation ideal filtration, and
we can form the maps
\begin{gather*}
\bar{g}_{2n} :F^{2n/2n+1}A= IA^{\,n}\big/IA^{\,n+1}\lra
IB^{\,n}\big/IB^{\,n+1} = F^{2n/2n+1}B,\\
\eo g:\eo A\lra \eo B.
\end{gather*}
Each $\bar{g}_{2n}$ is a map of graded modules, while $\eo A$ and $\eo B$ are
\emph{bigraded} algebras (by internal degree and filtration degree). It is easy to check
that $\eo g$ is surjective.  If $\bar{g}_{2n}(\bar{y})=0$, this means $g(y)\in
F^{2n+1}B$. Since $g$ is surjective, $y$ can be written as $y=a+b$, where $a\in
F^{2n+1}A$ and $b\in \ker g\cap F^{2n}A$. Since $\ker g=\image f$, $b=f(b')$ for some
$b'\in M$. Unfortunately, we cannot tell if $b'\in F^{2n}M$. If we refilter $M$ by
$\hat{F}^{2n}M=\inv{f}\bl F^{2n}A\br$, we obtain the exact sequence of bigraded algebras
$$\mathrm{gr}.^{\hat{F}}M\lra \eo A\lra \eo B.$$
In our case, $\gr^{\hat{F}} M\simeq\eo M \simeq M$, except that the generator of
$\gr^{\hat{F}} M$ has filtration degree equal to the filtration degree of $f(x)$ in $A$.
With the understanding that $M$ is regraded as such,
let us investigate the extension
\beqs
M \lra \eo A\lra \eo B
\enqs
of primitively generated Hopf algebras with associated Cartan-Eilenberg spectral sequence
\beq\label{yet another spectral sequence}
E_2^{s,t}\simeq H^s(\eo B)\ot H^t(M)\Lra H^{s+t}(\eo A).
\enq

We recall that a Hopf algebra $C$ is \emph{monogenic} if it
is generated by one element $x$.

\begin{proposition}\label{monogenic algebras}\cite{Milnor-Moore} A monogenic Hopf algebra
over $\f$ has the form
\be
\item\label{exterior monogenic} An exterior algebra $\Lambda[x]=\f[x]/(x^2)$, with
$\deg x$ odd for $p$ odd.
\item A truncated polynomial algebra $A=\f[x]\big/\bigl(x^{p^n}\bigr)$ as in \eqref{truncated
polynomials}.
\item A polynomial algebra $\f[x]$ with $\deg x$ even for $p$ odd.
\ee
\end{proposition}
In all cases, $x$ is primitive by reason of dimension. We have excluded the case $p=2$,
$n=1$ from $(2)$, as $\ff[x]/(x^2)\simeq \Lambda[x]$.
\begin{lemma} {\cite[Proposition 2.2]{Wilkerson}} For each monogenic Hopf algebra $C$ above,
\be
\item $H^*(C)\simeq\f[z]$, with $z\in H^{1,\deg x}$.
\item $H^*(C)\simeq\Lambda[z]\ot\f[e]$, with $z\in H^{1,\,\deg x}$, $e\in H^{2,\,p^n \! \deg x}$.
\item $H^*(C)\simeq\Lambda[z]$, with $z\in H^{1,\deg x}$.
\ee
\end{lemma}

\begin{theorem}[Borel, \cite{Milnor-Moore}] A locally finite, connected, primitively
generated, commutative Hopf algebra over $\f$ is
isomorphic to a tensor product of monogenic Hopf algebras.
\end{theorem}

\begin{theorem}\label{Hopf algebra main theorem} Let $M,A,B$ be as in \eqref{general extension}. If $B$ is
commutative, then $\eo A$ is semi-Koszul.
\end{theorem}
\begin{proof} The fiber $M$ has cohomology generated by $z\in H^1(M)$ and possibly
$e\in H^2(M)$. The element $e$ is related to $z$ by the Steenrod operation $\theta$,
where $\theta$ is one of $\bp$ or $\sqit{1}$. The spectral sequence \eqref{yet another
spectral sequence} will collapse if $d_3(e)=\theta\bl d_2(z)\br=0$. Since $B$ is
commutative, so is $\eo B$. By Borel's theorem, $\eo B\simeq \otimes_i B_i$, where the
$B_i$ are monogenic Hopf algebras.
Then $H^*(\eo B)\simeq H^*(B_1)\ot\dotsb\ot H^*(B_n)\simeq \Lambda[z_i]\ot\f[z_j']\ot
\f[e_k]$ with $z_i,\, z_j'\in H^1$, $e_k\in H^2$. We can show $\theta \equiv 0$ on
$H^2(\eo B)$ by using the same argument as in Theorem \ref{base is abelian}.
\end{proof}

\begin{remark}
A theorem of Wilkerson \cite{Wilkerson} states that a \emph{finite dimensional}, graded,
connected, cocommutative Hopf algebra $A\neq k$ always contains a nontrivial central
monogenic sub-Hopf algebra $C$. In the case that $C$ is a truncated polynomial algebra of
height $p^n$, we can always find a central primitive element $x\in A$ of height $p$. The
conditions of Theorem \ref{Hopf algebra main theorem} then hold.
\end{remark}


\section{Rank Two $p$-Groups}\label{rank 2}

In this section we show that the cohomology of the Lie algebra associated with a rank $2$
$p$-group for $p\geq 5$ is bigenerated.  A rank 2 $p$-group is a finite $p$-group whose
largest elementary abelian subgroup has dimension 2. In the case $p\geq 5$, Blackburn
\cite{Blackburn} has classified these as belonging to one of the following four families:
abelian, metacyclic, $C(r)$ and $G(r,e)$ of order $p^r$. A metacyclic $p$-group $P$ is an
extension of a cyclic group by a cyclic group. A presentation for $C(r)$ is given by
$$C(r)=\langle a,b,c \; \vert \; a^p=b^p
=c^{p^{r-2}}=1,\, (a,b)=c^{p^{r-3}}, \, c\in \mathcal{Z}(C)\rangle$$
where $r\geq 3$, $(a,b)=\inv{a}\inv{b}ab$, and $\mathcal{Z}(C)$ is the center of $C$. A
presentation for $G(r,e)$ is given by
$$G(r,e)=\langle f,g,h \; \vert \; f^p=g^{p^{r-2}}=h^p= \lp g,h\rp=1,\, \lp
f,h^{\sss{-1}}\rp=g^{ep^{r-3}},\, \lp f,g\rp=h\rangle$$ where $r\geq 4$,\, $e$ is 1 or a
quadratic non residue mod $p$.

These also define rank 2 $p$-groups for $p=3$. We include the families for $p=3$ in the
calculations to follow because they are the source of most of the interesting behavior.

We will be performing explicit calculations in the cobar complex and with Massey
products. See the Appendix for a review of these concepts.
\subsection{$G$ of type $C(r)$}\label{section 1}\mbox{}

Since $c^{p^{r-3}}$ is a central element of order $p$, we have the central
extension

$$\z\lra G\lra \z\times\z\times\zn{r-3}.$$
As the base is abelian, Theorem \ref{base is abelian} tells us that $\VL G$
is semi-Koszul.  We shall now explicitly calculate $H^*(\VL G)$.

\subsubsection{Case 1: $r\geq4$}

\begin{alignat*}{2}
\notag\Gamma^2 &= \bigl\langle a, b, c \bigr\rangle = G &\qquad \Gamma^{2/3}
&= \f\langle \overline{a},\overline{b},\overline{c} \rangle\\
\notag\Gamma^{2p} &= \langle c^p \rangle \simeq \mathbb{Z}/(p^{r-3})
&\qquad
\Gamma^{2p/2p+1} &= \f\langle \overline{c^p}\rangle\\
& \, \vdots & & \, \vdots \\
\notag\Gamma^{2p^{r-3}} &= \langle c^{p^{r-3}}\rangle \simeq
\z &\qquad \Gamma^{2p^{r-3}/2p^{r-3}+1} &=
\f\langle\overline{c^{p^{r-3}}}\rangle
 \\
\notag & & \qquad \Gamma^{i/i+1} &=0 \text{ for } i\neq 2 p^j,\, j\leq r-3.
\end{alignat*}

Since $a^{-1}b^{-1}ab= c^{p^{r-3}} \in \Gamma^{2p^{r-3}}$ falls below $\Gamma^4$, $\L G$
is an abelian restricted Lie algebra over $\f$.  Moreover, $\xi\bigl(
\overline{a}\bigr)=\xi\bigl(\ol{b}\bigr)=0,\, \xi\bigl(\overline{c^{p^i}}\bigr)=
\overline{c^{p^{i+1}}}$. Let $x=\ol{a},\, y=\ol{b},\, z_i=\overline{c^{p^{i-1}}}$ for
$1\leq i\leq r-2$. Then
\begin{align*}
\VL G&=\Tens\{x,y,z_i\} \big/
 \left(\begin{smallmatrix}
   x^p=y^p=0, \, z_i^p=z_{i+1}\\
   z_{r-2}^p=0, \,  x,y,z_i \text{ commute }
 \end{smallmatrix}\right) \\
&\simeq \f[x,y,z]\big/\bl x^p,y^p,z^{p^{r-2}}\br \quad \text{as Hopf algebras.}
\end{align*}
So $H^* ( \VL C(r) )\simeq \Lambda[x,y,z]\otimes\f[u,v,w] \up{where} x,y,z\in
\Ext^{1,2}$, and $u,v\in Ext^{2,2p}$ and $w\in Ext^{2,2p^{r-2}}$.  As usual, the
one-dimensional classes $x,y$ are related to the two dimensional classes $u,v$
 by the algebraic Bockstein $-\bp$, and $w=\langle z\rangle^{p^{r-2}}$.  See Appendix
\ref{calculations} for details.
\smallskip

\subsubsection{Case 2: $r=3$}\mbox{}

Let $G=C(3)=\bigr\langle a,b,c \, \vert
\, a^p=b^p =c^p=1,\, (a,b)=c\in \mathscr{Z}(G)\bigr\rangle$.
This case is more interesting, as the filtration is different.
The filtration degree of $c$ is 4, which makes $\L G$ non-abelian: $\bigl[ \ol{a},\ol{b}
\bigr]=\ol{c}\in\Gamma^{4/5}$. We still have  $\xi(\overline{a})=\xi\bigl(\ol{b}\bigr)=0$.
Thus $\VL G \simeq \Tens\{a,b,c\}/J$, where $J$ is the ideal generated by the relations
$\bigl\{ab-ba=c, \,ac=ca, \,bc=cb, \,a^p=b^p=c^p=0\bigr\}$. The central extension of groups:
\beq\label{group extension for C}
\begin{CD}
\Z/(p) @>>z\mapsto c> G @>>a,b\mapsto x,y> \Z/(p) \times \Z/(p)
\end{CD}
\enq
induces the central extension of algebras
\beq\label{algebra extension for C}
\begin{CD}
\f[z]\big/(z^p) @>>z\mapsto c> VLG @>>a,b\mapsto x,y> \f[x,y]\big/(x^p,y^p)
\end{CD}
\enq
which we abbreviate as $F \lra A \lra B$. This sequence will be exact as \emph{graded} algebras
if we set $\deg z=4$.  The induced Cartan-Eilenberg spectral sequence has
$E_2$ term
\begin{equation}\label{ss for C}
E_2^{*,*}\bl C(3)\br=H^*(B)\otimes H^*(F) \simeq
\Lambda[\sigma,\tau]\!\otimes\!\f[\mu,\nu]\, \otimes\, \Lambda[
\zeta]\!\otimes\!\f[\epsilon]
\end{equation}
$$\begin{aligned}[t]
\sigma, \tau &\in E_2^{1,0,2} \ &\mu, \nu&\in E^{2,0,2p}_2 \\
\zeta &\in E_2^{0,1,4}   &\epsilon &\in E_2^{0,2,4p}
\end{aligned}
\qquad
\begin{aligned}[t]
\mu =& -\bp(\sigma)=\langle\sigma\rangle^p, \quad \nu = -\bp(\tau)=\langle\tau\rangle^p\\
\epsilon =& -\bp(\zeta)=\langle\zeta\rangle^p
\end{aligned}$$
\smallskip

The third superscript is the internal degree. The elements $\sigma\otimes 1$ and
$\tau\otimes 1$ in $E^{1,0}_2$ are permanent cycles converging to nonzero classes
$\alpha$ and $\beta$ in $H^{1,2}(A)$. Put a reverse lexicographic order on $A:\, c>b>a$.
Choose the basis for $A$ induced by this ordering: $\mathscr{B}=\bigl\{c^i b^j a^k \
\vert \, 0\leq i,j,k\leq p-1\bigr\}$.  The nontrivial relations are generated by the
following:
\beqs
ab=ba+c, \quad ac=ca, \quad bc=cb, \quad a^p=b^p=c^p=0.
\enqs
Let $\mathscr{B}^*=\{x^*\vert x\in \mathscr{B}\}$ denote the dual basis, i.e.,
\beqs
\langle x^*, y\rangle =
\begin{cases}
1 &\quad \text{if } y=x,\\ 0&\quad \text{otherwise.}
\end{cases}
\enqs
Identify $\alpha$, $\beta$ with the classes $\bigl[a^*\bigr]$, $\bigl[b^*\bigr]$ in the
cobar complex $C^*(A)$. Now, $\delta\bigl[c^*\bigr]=\bigl[a^*\vert b^*\bigr]\,
\Rightarrow\, \alpha\beta=0 \up{in} E_\infty^{2,0}$.  This implies
\beq\label{d2 for C}
d_2(1 \otimes\zeta)=\lambda \cdot \sigma\otimes\tau
\enq
for some $\lambda\in\f^\times$,
because that is the last differential which can hit $E_2^{2,0}$.  We can ignore
the constant $\lambda$, all that matters is that it is nonzero.  Since $\sigma\ot\tau=
1\ot\sigma\cdot 1\ot\tau$ is a product, \eqref{proof that bp} implies that
the spectral sequence \eqref{ss for C} collapses at $E_3$. Moreover, it
can be easily seen that
\be
\item $d_2:\Tens (\mu,\nu)\otimes\zeta\epsilon^k\lra
\Ideal(\sigma\tau)\otimes\epsilon^k$ is never zero,
\item $d_2 \equiv 0$ on $\Ideal(\sigma,\tau)\ot\zeta\epsilon^k \subset E_2^{*,2k+1},$
\item $d_2\equiv 0$ on the $\epsilon^k$ lines.
\ee

The permanent classes
\begin{gather}
S=\{\sigma\ot 1, \, \tau\ot 1, \, \mu\ot 1, \, \nu \ot 1, \, \sigma \ot \zeta, \,
\tau \ot \zeta, \, 1 \ot \epsilon\}
\label{generators 1}
\intertext{in $E_\infty$ converge to elements}
T=\{\alpha,\beta\, u, v, x_1, x_2,e \} \label{generators 2}
\end{gather}
in $H^*(A)$, with $\alpha,\,\beta\in H^{1,2}(A),\: u,\,v\in H^{2,2p}(A),\:
x_1,\,x_2\in H^{2,6}(A)$ and $e\in H^{2,4p}(A)$.
Every element in $E_\infty$ can be expressed as a sum of products from
$S$, for example, $\sigma\mu^i\nu^j\ot\zeta\epsilon^k=
(\sigma \ot \zeta)(\mu\ot 1)^i (\nu\ot 1)^j (1\ot\epsilon)^k$. Therefore,

\begin{theorem}
$T$ generates $H^*(A)=H^*\bl\VL C(3)\br$ multiplicatively.
\end{theorem}

We shall look at the relations in $H^*(A)$ in the following section.

\subsection{Relations in $H^*\bl\VL C(3)\br$}\mbox{}

We will use a theorem of Ravenel \cite[Appendix 1.4]{Ravenel} on Massey products in
spectral sequences.  Refer to \ref{massey products} for the definition of Massey
products.  The formulation is somewhat complicated, but the main idea, due to Kahn
\cite{Kahn}, is that we do not want any higher ``cross-differentials'' to interfere with
$d_2$. Since $d_r=0$ for $r\geq 3$, this is not an issue.  Thus we can safely state that
if $a_i\in E_2^{*,*}$ converge to $u_i$ in $H^*(A)$, then the Massey product $\langle
a_1,\dots,a_n \rangle\in E_3$ (when defined) will converge to $\langle u_1,
\dots,u_n\rangle$, up to indeterminacy of course. We shall also make use of the following
``juggling theorem":
\begin{lemma} {\cite[A1.4.6]{Ravenel}}
If  $\langle u_1,u_2, u_{3}\rangle$ and  $\langle u_2,u_3, u_4\rangle$
are defined, then
\beqs
u_1 \langle u_2,u_3, u_4\rangle=\langle \ol{u}_1,\bar{u}_2, \ol{u}_{3}\rangle u_4,
\enqs
where $\ol{u}=(-1)^{\deg u+1}u$ and $\deg u$ means total degree of $u$.
\end{lemma}

\begin{proposition}\label{relations} In $H^*(\VL C(3))$, the
multiplicative relations are generated by the following:

\begin{enumerate}
\item $\alpha^2=\alpha\beta=\beta^2=0$.\label{eq0}
\item $\alpha x_2=-\beta x_1.$\label{eq1}
\item $\alpha x_1=
\begin{cases}
\beta u,\quad &p=3,\\
0,\quad &p>3.
\end{cases}\qquad
\beta x_2=
\begin{cases}
-\alpha v,\quad &p=3,\\
0,\quad &p>3.
\end{cases}$\label{eq2}
\item $x_1^2,x_1 x_2, x_2^2=\begin{cases}-ux_2, -uv, vx_1, \quad &p=3,\\
0,0,0,\quad &p>3.\end{cases}$\label{eq3}
\end{enumerate}

\noindent In addition, we have the following Massey products

\begin{enumerate}\addtocounter{enumi}{4}
\item $x_1 = \langle\alpha, \alpha, \beta\rangle=-1/2\langle \alpha,
\beta, \alpha\rangle=\langle\beta, \alpha, \alpha\rangle.\label{eq4}$
\item $x_2=-\langle\alpha,\beta,\beta\rangle=1/2\langle\beta,\alpha,\beta\rangle
=-\langle\beta,\beta,\alpha\rangle.\label{eq5}$
\item $\langle\alpha\rangle^p=-\bp(\alpha)=u, \quad \langle\beta\rangle^p=-\bp(\beta)=v.
\label{eq6}$
\end{enumerate}
\end{proposition}

We will prove the last three first.

\noindent\textit{Proof of} \ref{relations}.\ref{eq4},  \ref{relations}.\ref{eq5}.

Using the notation of \S  \ref{massey products}, let $(U,\delta)=\bl E_2^{**}, d_2\br$.
The elements $a_1=a_2=\sigma\ot 1, \, a_3=\tau\ot 1$ are in $E_2^{1,0}$, and have total
degree 1. Therefore $\ol{a_i}=a_i$. Now $a_1 a_2=\sigma^2\ot 1=0,$ and $a_2 a_3
=\sigma\ot\tau = d_2(1\ot\zeta)$, so we can take $u_1=0$,\, $u_2=1\ot\zeta$. Then
$\langle a_1,a_2,a_3 \rangle=\sigma\ot 1\cdot 1\ot\zeta + 0\cdot\tau\ot 1 =
\sigma\ot\tau$, which converges to $x_1$.  The indeterminacy is $\alpha H^1 + H^1
\beta$, which is 0 since $H^1(A)=\f\langle\alpha,\beta\rangle$ and $\alpha^2=\alpha\beta=\beta^2=0$.

$\langle\beta, \alpha, \alpha\rangle=\tau\ot 1\cdot 0 +(-1\ot\zeta)\cdot\sigma\ot 1
=\sigma\ot\zeta$.
The minus sign comes from the fact that $\beta\alpha=-\alpha\beta$.  It disappears
because $1\ot\zeta$ and $\sigma\ot 1$ anticommute.

$\langle \alpha,\beta, \alpha\rangle=\sigma\ot 1\cdot(-1\ot\zeta)+1\ot\zeta\cdot
\sigma\ot 1=-2\,\sigma\ot\zeta$.

Now let $U=C^{**}(A)$, the cobar complex of $A$.  $\alpha$ is represented by $a_1=a_2=[a^*]$,
$\beta$ by $a_3=[b^*]$.  Let us calculate $\langle\alpha, \alpha, \beta\rangle$ in $U$.
$u_1=[{a^2}^*],\, u_2 =[c^*]$, and $\langle a_1,a_2,a_3 \rangle=[a^*\vert c^*]
+[{a^2}^*\vert b^*]$.  We now have a representative for $x_1$ in $C^*(A)$.

Other representatives are
\begin{subequations}\label{Massey for x1}
\begin{align}
x_1&=-\recip{2}\langle \alpha,\beta,\alpha\rangle=\frac{1}{2}\bl [a^*\vert c^*]-
[c^*\vert a^*]-[a^*\vert {ba}^*]\br\\
x_1&=\langle\beta,\alpha,\alpha\rangle=[b^*\vert{a^2}^* ] -[c^* \vert a^* ]+
[{ba}^* \vert a^* ]
\end{align}
\end{subequations}
The procedure for $x_2$ is similar.  Some representatives are
\begin{subequations}\label{Massey for x2}
\begin{align}
x_2&=-\langle\alpha,\beta,\beta\rangle=-[a^*\vert{b^2}^*]-[c^*\vert b^*]\\
x_2&=-\langle\beta,\beta,\alpha\rangle=-[b^*\vert{ba}^*]+[b^*\vert c^*]-[{b^2}^*\vert a^*]\\
x_2&=\recip{2}\langle\beta,\alpha,\beta\rangle=\recip{2}\bl[b^*\vert c^*]+[{ba}^*\vert b^*]
-[c^*\vert b^*]\br.
\end{align}
\end{subequations}\qed

\noindent\textit{Proof of} \ref{relations}.\ref{eq6}.

Follows from naturality of the algebraic Steenrod operations.

\noindent\textit{Proof of} \ref{relations}.\ref{eq0}.

Follows from relations in $E_\infty$.

\noindent\textit{Proof of} \ref{relations}.\ref{eq1}.

$\alpha x_2=-\alpha\langle\alpha,\beta,\beta\rangle=-\langle\ol{\alpha},\ol{\alpha},\ol{\beta}\rangle
\beta=-x_1 \beta =-\beta x_1$

\noindent\textit{Proof of} \ref{relations}.\ref{eq2}.

$\alpha x_1=\alpha\langle\alpha,\alpha,\beta\rangle=\langle\ol{\alpha},\ol{\alpha},
\ol{\alpha}\rangle\beta=\langle\alpha,\alpha,\alpha\rangle\beta=u\beta$
for $p=3$. For larger primes, $\langle\alpha\rangle^3=0$, since for
$\langle\alpha\rangle^p$ to be defined, we need all lower symmetric Massey
products to vanish. Similarly, $\beta
x_2=-\beta\langle\beta,\beta,\alpha\rangle=-\langle\beta,\beta,\beta\rangle
\alpha=-v\alpha$ for $p=3$, while it vanishes for larger primes. In either
situation, the products fall to the lower filtration level
$F^3H^3(A)=E_\infty^{3,0}$.\qed

\noindent\textit{Proof of} \ref{relations}.\ref{eq3}.

Since ${(\sigma\ot\zeta)}^2 =-\sigma^2\ot\zeta^2=0$ in $E_\infty^{2,2}\simeq
F^{2/3}H^4(A)=F^2 H^4(A)\big/F^3 H^4(A)$, we deduce that $x_1^2\in F^3 H^4(A)$.  Since
$F^3 H^4(A)$ is spanned by $\bigl\{u x_1, ux_2, vx_1, vx_2, u^2,$ $uv, v^2\bigr\}$, we
can express $x_1^2$ as a linear combination
\begin{equation}\label{expr}
x_1^2=c_1 u x_1+c_2 ux_2+c_3 vx_1+c_4 vx_2+c_5 u^2+c_6 uv +c_7 v^2, \quad
c_i\in \f.
\end{equation}
\noindent The element $x_1^2$ has internal degree $12$, while the elements on the right
hand side of (\ref{expr}) have internal degree larger than or equal to $2p+6$.  For
$p>3$, we are forced to conclude that $x_1^2=0$.  As usual, the case $p=3$ requires
special attention.

Let $p=3$. Multiplying both sides of (\ref{expr}) by $\alpha$ gives us $\alpha
x_1^2=\beta ux_1$ on the left, and $c_1\beta u^2-c_2\beta u x_1 +c_3\beta uv-c_4\beta
vx_1+c_5\alpha u^2+ c_6\alpha uv+c_7\alpha v^2$ on the right. The summands are linearly
independent in $H^5(A)$ because their representatives in $E_\infty$ are linearly
independent. Therefore, $c_2=-1$, and all the other $c_i$ are zero. Starting with linear
combinations in $H^4(A)$ for $x_1 x_2$ and $x_2^2$, then multiplying by $\alpha$, will
yield the other two identities.\qed

These relations all come about from the fact that $\sigma^2=\tau^2=\zeta^2=0,\,
\sigma\tau=-\tau\sigma$ in $E_\infty$.  Since there are no other identities of this type
in $E_\infty$, there are no new relations in $H^*(A)$.


\subsection{$G$ of type $G(r,e)$}\mbox{}

We next turn our attention to groups $G$ of type
$$G(r,e)=\langle f,g,h \; \vert \; f^p=g^{p^{r-2}}=h^p= \lp g,h\rp=1,\, \lp
f,h^{\sss{-1}}\rp=g^{ep^{r-3}},\, \lp f,g\rp=h\rangle.$$
 We will soon see that for our purposes, we can take $e=1$.
\begin{lemma}\label{center}
 $g^p\in\mathcal{Z}(G).$
\end{lemma}
\begin{proof} Since $(f,g)=h$, we have $\inv{f}\inv{g}f=h\inv{g}$, and
\beq
(f,g^{p})=f^{-1}g^{-p}f g^{p}={\lp f^{-1}g^{-1}f\rp}^{p} g^p =\bl hg^{-1}\br^p
g^p=h^p g^{-p} g^p=1
\enq
The second to last equality follows because $(g,h)=1$.
\end{proof}

We will investigate the exact sequence
\begin{equation}\label{yae}
\begin{CD}
\Z/(p) \langle d\rangle @>>d\mapsto g^{p^{r-3}}> G(r) @>>>
Q(r-1)
\end{CD}
\end{equation}
This is a central extension by Lemma \ref{center}. As the base is not abelian, we cannot
use Theorem \ref{base is abelian}.
The cokernel $Q$ has generators $\tilde{f},\tilde{g},\tilde{h}$ and relations:
\begin{equation*}
\tilde{f}^{\,p}=\tilde{h}^{\,p}=\bigm( \tilde{g},\tilde{h} \bigm)=1,\,
 \bigm( \tilde{f},\tilde{h}^{\,\scriptstyle{-1}} \bigm)=\tilde{g}^{\,ep^{r-3}}=1,\;
 \bigm( \tilde{f},\tilde{g} \bigm)=\tilde{h}.
\end{equation*}
Therefore $Q(r-1)\simeq\langle a,b,c \; \vert \; a^p=b^{p^{r-3}}=c^p=1,\,
(a,b)=c\in \mathscr{Z}(Q)\rangle,$ where $a=\tilde{f}$, $b=\tilde{g}$,
$c=\tilde{h}$. For $r=4$,  $Q(r-1)=Q(3)\simeq C(3)$. We will focus on this case
presently. The case $r>4$ is taken up in section \ref{case 2}.

\subsubsection{Case 1: $r=4$}\label{counterexample}\mbox{}

For the extension \eqref{yae}, we have
\begin{equation}\label{G4 extension}
\begin{CD}
  \z @>>d\mapsto g^p>  G(4) @>>f,g,h\mapsto a,b,c>   Q=C(3)
\end{CD}
\end{equation}
with $f,g,h,d$ in filtration degree $2,2,4,2p$ respectively. The Lie algebra $\L G$  is
generated by $\left\{ x=\ol{f}, y_1=\ol{g}, y_2=\overline{g^p}, z=\ol{h}\right\}$. Since
$2+4=2p\Leftrightarrow p=3$, $\lb\,\ol{f},\ol{h}\,\rb=0$ for $p>3$. Since
$\overline{\inv{h}}=-\ol{h}$,
$\lb\,\ol{f},\ol{h}\,\rb=-\lb\,\ol{f},\overline{\inv{h}}\,\rb=-\overline{g^{ep}}
=-e\overline{g^p}$ for $p=3$. Summarizing, the non-trivial structure of $\L G$ is as
follows
\begin{equation*}
\xi(y_1)=y_2,\quad [x,y_1]=z,\quad [x,z]=
\begin{cases}
-ey_2 \quad & p=3\\
 0         \quad & p>3.
\end{cases}
\end{equation*}
For $p=3$, $e$ can only equal 1, so we see that $e$ doesn't really matter.

Taking the usual quotient of $\Tens\{x,y_1,y_2,z\}$ by the ideal $J$ of relations as in
Definition \ref{VL},
\begin{equation}\label{VLG for G4}
\VL G \simeq
\begin{cases}
\Tens\{f,g,h\}\big/
   \left(\begin{smallmatrix}
f^p=g^{p^2}=h^p=0, \, fg=h+gf \\
    fg^p=g^p f,\, fh=hf,\, gh=hg\end{smallmatrix}\right) &\quad p>3,\\
\vspace{-3mm}\\
\Tens\{f,g,h\}\big/
   \left(\begin{smallmatrix} f^p=g^{p^2}=h^p=0,\, fg=h+gf\\
    fg^p=g^p f,\, fh=-g^p+hf,\, gh=hg\end{smallmatrix}\right)   &\quad p=3.
\end{cases}
\end{equation}
under the mapping $x\mapsto f,\, y_1\mapsto g,\, y_2\mapsto g^p,\, z\mapsto h$.

The central extension of Hopf algebras:
\begin{equation}\label{algebra extension for G4}
\begin{CD}
\f[d]\big/(d^p) @>i>d\mapsto g^p> \VL G @>\pi>f,g,h\mapsto a,b,c> \VL Q=\VL C(3).
\end{CD}
\end{equation}
denoted by $F\stackrel{i}{\lra} A\stackrel{\pi}{\lra} B$ gives us a spectral sequence
\begin{equation}\label{ss for G4}
E_2^{*,*}\bl G(4)\br\simeq H^*(B) \ot \left(\Lambda[z]\ot\f[w]\right) \Rightarrow H^*(A).
\end{equation}
where $z\in E_2^{0,1,2p}$ and $w\in E_2^{0,2,2p^2}$. We shall suppress tensor notation in
spectral sequences (e.g $1\ot z=z$) from now on.
\begin{proposition}\label{d2 of z} In the spectral sequence \eqref{ss for G4}
 $$d_2(z)=
\begin{cases}
v-x_1, \quad &p=3\\
v,     \quad &p>3.
\end{cases}$$
\end{proposition}
\begin{proof}
Put an ordering $h>g>f$ on $A=\VL G$. A basis for $A$ consists of $\{h^i g^j f^k\vert
0\leq i,k\leq p-1 ,\, 0\leq j\leq p^2-1\}$. Using the relations \eqref{VLG for G4} for
$A$, we see that in the cobar complex $C^*(A)$
\begin{equation}
\label{diff}
\delta [{g^p}^*]=\begin{cases}
\sum_{i=1}^{p-1}[{g^i}^*\vert {g^{p-i}}^*] -[f^*\vert h^*]-[{f^2}^*\vert g^*]&\quad p=3,\\
\sum_{i=1}^{p-1}[{g^i}^*\vert {g^{p-i}}^*]&\quad p>3.
\end{cases}
\end{equation}
Under the edge map $\pi^*: H^*(B)\lra H^*(A)$, $v\mapsto
\cl\bl\sum_{i=1}^{p-1}[{g^i}^*\vert {g^{p-i}}^*]\br$ and $x_1\mapsto
\cl\bl[f^*\vert h^*]+[{f^2}^*\vert g^*]\br$. The only differential that
can cause the right hand side in \eqref{diff} to be zero is $d_2(z)$.
\end{proof}
\begin{proposition}  \label{bp}
$\bp(v)=\bp(x_1)=0$ in $H^3(B)$.
\end{proposition}
\begin{proof}
$v, x_1 \in H^{2,2p}(B)$. The image of $\bp : H^{2,2p}(B)\lra H^{3,2p^2}$ is zero because the internal degrees
of the elements of $H^3(B)$ are all less than or equal to $2+4p$, which is less
than $2p^2$.
\end{proof}
\begin{proposition}\label{zero divisor}
$d_2(z)$ is not a zero divisor in $H^*(B)$.
\end{proposition}
\begin{proof} Assume $d_2(z)y=0$ for some nonzero element $y\in H^{s+t}(B)$ of filtration
degree $s$.  Project $y$ to $\bar{y}\in F^{s/s+1}H^{s+t}(B)\simeq E_\infty^{s,t}\bl
C(3)\br$.

For $p>3$, the equation $vy=0$ projects to $ \nu\bar{y}=0$. Since $\nu$ is not a zero
divisor in $E_\infty\left( C(3)\right)$, $\bar{y}=0$.  This means $y\in F^{s+1}H^*(B)$,
which is a contradiction.

For $p=3$, we have $(v-x_1)y=0$, or $vy=x_1y$. This projects to $\nu\bar{y}
=\sigma\zeta\bar{y}$. Note that $\nu\bar{y}\in E_\infty^{s+2,t}$, while $\sigma
\zeta\bar{y}\in E_\infty^{s+1,t+1}$.  The two cannot possibly be equal, unless of course
they are both zero.  But $\nu$ is not a zero divisor, and the same contradiction ensues.
\end{proof}

We conclude
\begin{proposition}\label{result}
The spectral sequence \eqref{ss for G4} collapses at $E_3$, and
$$H^*(A)\simeq
\begin{cases}
 H^*(B)\big/(v-x_1)\ot\f[w],\quad &p=3,\\
 H^*(B)\big/(v)\ot\f[w],\quad &p>3.
\end{cases}$$
\end{proposition}

\begin{theorem}\label{cohomology of G4}
$H^*(\VL G(4))$ is generated by $\{\alpha, \beta, u, x_1, x_2, e, w\}$, with
$\alpha,\beta\in H^{1,2}$,
$x_1, x_2\in H^{2,6}$, $u\in H^{2,2p}$,  $e\in H^{2,4p}$, and $w\in
H^{2,2p^2}$.
The relations are generated by
\begin{enumerate}
\item $\alpha^2=\alpha\beta=\beta^2=0$. 
\item $\alpha x_2=-\beta x_1.$ 
\item $\alpha x_1=
\begin{cases}
\beta u,\quad &p=3,\\
0,\quad &p>3.
\end{cases}\qquad
\beta x_2=
\begin{cases}
-\beta u,\quad &p=3,\\
0,\quad &p>3.
\end{cases}$  
\item $x_1^2,x_1 x_2, x_2^2=\begin{cases}-ux_2, -ux_1, -ux_2, \quad &p=3,\\
0,0,0,\quad &p>3.\end{cases}$
\end{enumerate}
In particular, $\VL G(4)$ is semi-Koszul. The Massey products of Proposition \ref{relations}
still hold, except now
$\langle\beta\rangle^p=x_1$ for $p=3$, and
$\langle\beta\rangle^p=0,\ \langle\beta\rangle^{p^2}=w$ for $p>3$.
\end{theorem}
\begin{proof} This is simply a matter of setting $v=0$ or $v=x_1$ in the relations
for $H^*\bl \VL C(3)\br$ (Proposition \ref{relations}). The class $w$ is
formally the coboundary of $\bigl[{g^{p^2}}^*\bigr]$. For $p>3$, this is
equal to the representative of $\langle\beta\rangle^{p^2}$ in the cobar
complex.
\end{proof}

\begin{remark}
Notice that $x_2^2=- u x_2$ implies $x_2+u$ is a zero divisor, with $x_2(x_2+u)=0$.
If we had a central extension of 3-groups
$\Z/(3) \lra K \lra G(4)$
and the induced spectral sequence \eqref{algebra spectral sequence} had
$d_2(z)=x_2+u$
then $x_2\ot z\in E_\infty^{2,1}$ would be an indecomposable element of dimension 3.  In this
case, $\VL K$ would not be semi-Koszul.
\end{remark}


\subsubsection{Case 2: $r>4$}\label{case 2}\mbox

Now let $G=G(r,e),\, r>4,\, p\geq 3$. The situation is
\begin{gather*}
\begin{CD}
\Z/(p)  \langle d\rangle @>>d\mapsto g^{p^{r-3}}> G(r) @>>f,g,h\mapsto a,b,c>
Q(r-1)
\end{CD}\\
\text{where }\, Q(r-1)\simeq\langle a,b,c \; \vert \; a^p=b^{p^{r-3}}=c^p=1,\,\lb a,b\rb=c\in \mathscr{Z}(Q)\rangle
\end{gather*}
The filtration degrees of $f,g,h,d$ are $2,2,4,2p^{r-3}$ respectively.

Let us first deal with finding $H^*(\VL Q(r-1))$. There is a central extension for $Q$
\beqs
\begin{CD}
\z\langle z\rangle @>>z\mapsto c> Q\langle a,b,c\rangle @>>a,b\mapsto x,y> \z
\times \zn{r-3}\langle x, y\rangle
\end{CD}
\enqs
which is very similar to \eqref{group extension for C}. The induced algebra
extension
\beqs
\begin{CD}
\f[z]\big/(z^p) @>>> \VL Q(r-1) @>>> \f[x,y]\big/\bl x^p,y^{p^{r-3}}\br, \quad \deg z=4
\end{CD}
\enqs
gives us the spectral sequence
$$E_2^{**}\bl Q(r-1)\br
\simeq \Lambda[\sigma,\tau]\!\otimes\!\f[\mu,\nu]\, \otimes\, \Lambda[
\zeta]\!\otimes\!\f[\epsilon]$$ with $\sigma,\tau,\mu,\nu,\zeta,\epsilon$ as in \eqref{ss
for C} except $\nu$ has internal degree $2p^{r-3}$ and $\nu=\langle
\tau\rangle^{p^{r-3}}$. Just as in \eqref{d2 for C},\, $d_2(\zeta)=\sigma\tau$ and the
spectral sequence collapses at $E_3$.
\begin{theorem}
$H^*\bl\VL Q(r-1)\br$ is generated by $\{\alpha,\beta,x_1,x_2,u,v,e\}$, with
$\alpha,\beta\in H^{1,2}$, $x_1, x_2\in H^{2,6}$, $u\in H^{2,2p}$, $v\in H^{2,2p^{r-3}}$
and $e\in H^{2,4p}$. Relations are generated by
\begin{enumerate}
\item $\alpha^2=\alpha\beta=\beta^2=0$.
\item $\alpha x_2=-\beta x_1.$
\item $\alpha x_1=
\begin{cases}
\beta u,\quad &p=3,\\
0,\quad &p>3.
\end{cases}\qquad
\beta x_2=0.$
\item $x_1^2,x_1 x_2, x_2^2=\begin{cases}-ux_2,0, 0, \quad &p=3,\\
0,0,0,\quad &p>3.\end{cases}$
\end{enumerate}
The Massey products of Proposition \ref{relations}
still hold, except $\langle\beta\rangle^{p^{r-3}}=v$ for all $p$.
%
\end{theorem}
\begin{proof}
The proofs of  \eqref{eq0} and \eqref{eq1} are the same as in \ref{relations}. For the
second half of \eqref{eq2}, $\beta x_2=-\beta\langle
\beta,\beta,\alpha\rangle=-\langle\beta,\beta,\beta\rangle\alpha=0$ even for
$p=3$.
In the proof of \eqref{eq3} for $p=3$, we are reduced to solving the equation
\beq\label{expr2}
x_ix_j=c_1 ux_1+c_2ux_2+c_5u^2, \quad  c_k\in \f,
\enq
by reason of internal dimension.  Multiplying both sides by $\alpha$ and equating coefficients yield
the desired identities.  For example, $\alpha x_2^2=-\beta x_1 x_2=0$, while the right side of
\eqref{expr2} becomes $c_1\beta u^2-c_2\beta u x_1+c_5 \alpha u^2$.
\end{proof}

We now turn to finding $H^*(\VL G(r))$. Notice that since
$(f,h^{-1})=g^{p^{r-3}}$ falls below $\Gamma^6 G$, $\lb\ol{f},\ol{h}\rb=0$ in
$\L G$ even for $p=3$.
Thus \eqref{VLG for G4} becomes
\beqs
\VL G(r)\simeq \Tens\{f,g,h\}\Big/\binom{f^p=g^{p^{r-2}}=h^p=0,\,fg=h+gf}
{ fg^p=g^pf, \, fh=hf,\, gh=hg }.
\enqs
As in \eqref{algebra extension for G4}, setting $\deg d=2p^{r-3}$ gives us the graded
 algebra extension
\beqs
\begin{CD}
\f[d]/(d^p) @>i>> \VL G(r) @>\pi>> \VL Q(r-1)
\end{CD}\enqs
which induces the analogue of \eqref{ss for G4}
\beq\label{ss for G(r)}
E_2^{*,*}\bl G(r)\br\simeq H^*\bl\VL Q(r-1)\br \ot \Lambda[z]\ot\f[w].
\enq
\begin{proposition} In the spectral sequence \eqref{ss for G(r)},
$d_2(z)=v$.
\end{proposition}
\begin{proof}
In the cobar complex $C^*\bl\VL G(r)\br$, $\delta \bigl[{g^{p^{r-3}}}^*\bigr]=\sum_{i=1}^{p^{r-3}-1}
\bigl[{g^i}^*\vert {g^{p^{r-3}-i}}^*\bigr]$, which is the representative of
$\pi^*(v)$.
Thus $d_2(z)=v$ for all $p$.
\end{proof}
Proceeding in the same way as in Propositions \ref{bp}, \ref{zero divisor} and
\ref{result}, we obtain
\begin{proposition}
$\bp(v)=0$, and $v$ is not a zero divisor in $H^*\bl\VL Q(r-1)\br$. Thus the spectral sequence
\eqref{ss for G(r)} collapses at $E_3$ and
\beqs
H^*\bl \VL G(r)\br\simeq E_\infty \simeq H^*\bl\VL Q(r-1)\br\big/(v)\otimes
\f[w]
\enqs
\end{proposition}
Finally,
\begin{theorem}
$H^*\bl\VL G(r)\br$ is generated by $\alpha, \beta, u, x_1, x_2, e, w$. The bidegrees
of the generators are as in Theorem \eqref{cohomology of G4}, except $w\in
H^{2,2p^{r-2}}$.
Relations are as in Theorem \ref{cohomology of G4} except now $\beta x_2=x_1
x_2=x_2^2=0$, $\langle\beta\rangle^{p^{r-2}}=w$ for all $p$.
In particular, $H^*\bl \VL G(r)\br$ is semi-Koszul.
\end{theorem}


\subsection{$G$ Metacyclic}

Let $G$ be a metacyclic $p$-group. That is, there is an extension
$$
\begin{CD}
\zn{m}@>>> G @>>> \zn{n}
\end{CD}
$$
For $p$ odd,  any  metacyclic $p$-group $G$ has presentation
\beqs
G =P(m,n,q,l)=\langle x,y \; \vert \; x^{p^m}=1,\, y^{p^n}=x^{p^q},\,
(x,y)=x^{p^l}\rangle
\enqs
for integers $m,n,l,q\geq 1$ satisfying $(p^l+1)^{p^n}\equiv 1 \mod p^m$ and
$(p^l+1)p^q\equiv p^q\mod p^m$\, \cite{Dietz-odd}.  We can get all groups up to
isomorphism by requiring $l\leq m$, $q\leq m$, $n+l\geq m$ and $q+l\geq m$. Note that
$y^{p^{m+n-q}}=1$.

It is clear that $\L G$ is generated by $x_i=\overline{x^{p^{i-1}}}\in
\Gamma^{2p^{i-1}/2p^{i-1}+1},\, 1\leq i\leq m$, and $y_j=\overline{y^{p^{j-1}}}\in
\Gamma^{2p^{j-1}/2p^{j-1}+1}$,
$1\leq j\leq n+m-q$, some of which may be zero.
\begin{lemma}\label{metacyclic lemma 1} For all $i,j$, \,$[x_i,y_j]=0$.
\end{lemma}
\begin{proof}
$\inv{y}xy=x^{p^l+1}\, \Ra\, y^{-p^j}x^{p^i}y^{p^j}=x^{p^i(p^l+1)^{p^j}}\, \Ra\, \bl
x^{p^i},y^{p^j}\br=x^{p^i\left( (p^l+1)^{p^j}-1\right)}.$ The exponent is $
p^i\bl(p^l+1)^{p^j}-1\br=p^i\bl p^{lp^j}+\dots +p^j p^l+1-1\br=p^{i+j+l}a,\quad a\in \Z.
$ Therefore, $\bl x^{p^i},y^{p^j}\br=(x^a)^{p^{i+j+l}}\in \Gamma^{2p^{i+j+l}}\subsetneq
\Gamma^{2(p^i+p^j)}$, even for $i=j=0$.
\end{proof}

 We can break the analysis of $\L G$ down into three cases:

\noindent\textbf{Case 1:}\, $n<q$.

$\xi(y_j)=y_{j+1}$ for $j< n$. Since $y^{p^n}=x^{p^q}$ falls down to a lower
filtration level, $\xi(y_{n})=\overline{y^{p^n}}=0$, and there are no $y_j$ for $j> n$. The other
generators
$\{x_i\}_{1\leq i\leq m}$ have restriction $\xi(x_i)=x_{i+1}$ for $i<m$ and
$\xi(x_m)=0$. Then
\beqs
\VL G\simeq \f[x,y]\big/\bigl(x^{p^m},y^{p^n}\bigr) \quad \deg x=\deg y=2.
\enqs

\noindent\textbf{Case 2:}\, $n>q$.

By the same reasoning as above, t here are no $x_i$ for $i>q$. $\xi(x_i)=x_{i+1}$ for $i<
q$, $\xi(x_{q})=0$. $\xi(y_j)=y_{j+1}$ for $j< n+m-q$, and $\xi(y_{n+m-q})=0$. So
\beqs
\VL G\simeq \f[x,y]\big/\bl x^{p^q},y^{p^{n+m-q}} \br \quad \deg x=\deg y=2.
\enqs

\noindent\textbf{Case 3:}\, $n=q$.

In this case, $x^{p^q}=y^{p^n}$ in $\Gamma^{2p^q}=\Gamma^{2p^n}$. Let
$z=y x^{p^{m-n}-1}$.
Then
$(x,z)=\inv{x}x^{1-p^{m-n}}\inv{y}x y x^{p^{m-n}-1}=x^{p^l}=(x,y)$.
\begin{claim} $z^{p^n}$=1.
\end{claim}
\begin{proof}
It can be shown by induction that $(yx^a)^{b}=y^b x^{a[(p^l+1)^b-1]/p^l}$. For
$a=p^{m-n}-1,\, b=p^n$, the exponent of $x$ is
$$\bl p^{m-n}-1 \br\bl (p^l+1)^{p^n}-1\br/p^l=\bl p^{m-n}-1 \br
p^{n+l}c/p^l =( p^{m}-p^{n})c,$$ where
$c=1+\recip{p^n}\dbinom{p^n}{2}p^l+\recip{p^n}\dbinom{p^n}{3}p^{2l}+\dotsb =1+p^l d$, for
some $d\in\Z$. Then $z^{p^n}=y^{p^n} x^{( p^{m}-p^{n})c}=x^{p^n (-c+1)}=x^{p^{n+l}d}=1$
since $n+l\geq m$.
\end{proof}

So $G\simeq \langle x,z\vert x^{p^m}=1,\, z^{p^n}=1,\, (x,z)=x^{p^l}\rangle$.
Generators for
$\L G$ are $\{x_i\}_{1\leq i\leq m}$ and $\{z_j\}_{1\leq j\leq n}$, and
\beqs
\VL G\simeq \f[x,z]\big/\bl x^{p^m},z^{p^{n}}\br \quad \deg x=\deg z=2.
\enqs
Note that the transformation $y\mapsto z$ in $G$ corresponds to the map $y_1\mapsto
z_1=y_1-x_1$ in $\L G$.

\begin{theorem} Let $G=P(m,n,q,l)$, then
\begin{gather*}
H^*(\VL G)\simeq \Lambda[x,y]\ot\f[u,v],\ x,y\in H^{(1,2)},\,  u=\langle x\rangle^a\in
H^{(2, 2a)},\,  v=\langle y\rangle^b\in H^{(2,2b)}
\end{gather*}
where $a=p^m,\, b=p^n$ if $q\geq n$ and $a=p^q,\, b=p^{m+n-q}$ if $q<n$.
\end{theorem}

We have finally proven
\begin{theorem} Let $p\geq 5$, and let $G$ be a finite $p$-group of rank 2. Then
$\VL G$ is semi-Koszul.
\end{theorem}

\appendix
\section{}
\subsection{The Cobar Complex}\label{cobar complex}

We summarize the construction of the \emph{cobar complex} of $A$, which can be used to
calculate the cohomology of $A$.  Refer to Priddy \cite[Section 1]{Priddy} for more
details.

As in Section \ref{algebraic preliminaries}, let $A$ be an algebra with augmentation
ideal $I$. The cobar complex of $A$, denoted by $C^*(A)$, is given by $C^s(A)=(I^*)^{\ot
s}$. A typical element is written $\bigl[\alpha_1\vert\dotsb\vert\alpha_s\bigr]$ and has
bidegree $(s,t)$, where $t=\sum_{i=1}^s \deg\alpha_i$ is the internal degree inherited
from $A$.  Let $C^{s,t}$ be spanned by elements of bidegree $(s,t)$. There is a
differential $\delta : C^{s,t}(A)\longrightarrow C^{s+1,t}(A)$, which may be computed as
follows. Let $\mu^*$ be the composite $A^* \stackrel{\mu^*_A}{\lra} \bigl(A\otimes
A\bigr)^* \stackrel{\theta^{-1}}{\lra} A^*\otimes A^*$, where $\mu^*_A$ is dual to the
multiplication map $\mu_A \text{ of } A$, and $\theta$ is the isomorphism defined more
generally by
\begin{equation*}
\theta :M^*\ot N^*\lra \bl M\ot N\br^*, \quad \theta(f\ot g)(m\ot n)=(-1)^{\deg g\,\deg
m} f(m)g(n)
\end{equation*}
for $A$-modules $M$ and $N$. If $\mu^*(\alpha)=\sum_r \alpha^{'}_r \otimes
\alpha^{''}_r$, then $\delta$ is given by
\begin{equation}\label{delta}
\delta \bigl( \bigl[ \alpha_1 \vert \dotsb \vert \alpha_s \bigr] \bigr) = -\sum_{1\leq i
\leq s; r} (-1)^{e_{i,r}} \bigl[ \alpha_1 \vert \cdots \vert \alpha^{'}_{i,r} \vert
\alpha^{''}_{i,r} \vert \cdots \vert\alpha_s \bigr],
\end{equation}
where $e_{i,r}$ is the total degree of $\bigl[\alpha_1 \vert \cdots \vert
\alpha^{'}_{i,r} \bigr]$. The cohomology groups of $A$ are then $H^{s,t}(A)=H^s\bl
C^{*,t}(A),\delta\br$.

$C^*(A)$ has the structure of a differential graded algebra with product
\begin{equation*}
[\alpha_1\vert \dotsb\vert \alpha_s]\cup [\beta_1\vert \dotsb\vert \beta_{s'}]=
[\alpha_1\vert \dotsb\vert \alpha_s \vert \beta_1 \vert \dotsb\vert \beta_{s'}].
\end{equation*}
The induced product in cohomology coincides with the Yoneda Ext product as given in
\cite[2.2]{Adams}.

\subsection{Massey Products}\label{massey products}
The following material on Massey products can be found in \cite[A1.4]{Ravenel} and
\cite{Kraines}. We recall the definition of the 3-fold Massey product . Let $(U,\delta)$
be a differential bigraded algebra. Let $\alpha_i \in H^{s_i,t_i}(U)$ be represented by
$a_i$, and assume $\alpha_1\alpha_2= \alpha_2\alpha_3=0$. Let $s=s_1+s_2 +s_3, \, t=t_1+
t_2 + t_3$. For any element $a\in U$, let $\bar{a}=(-1)^{\deg(a)+1}a$ where, as usual,
degree means total degree.  Let $\delta\bl u_i\br=\bar{a}_i a_{i+1}$. Then the Massey
product $\langle\alpha_1,\alpha_2,\alpha_3\rangle\in H^{s-1,t}(U)$ is the class
represented by the cocycle $\bar{a}_1 u_2+\ol{u}_1a_3$. This 3-fold product is not well
defined because the choices made in its construction are not unique. The choices of $a_i$
don't matter, but the $u_i$ could each be altered by adding a cocycle $x_i\in
U^{s_i+s_{i+1}-1, t_i+t_{i+1}}$. This would alter
$\langle\alpha_1,\alpha_2,\alpha_3\rangle$ by an element of the form $\alpha_1 \llbracket
x_2\rrbracket+\llbracket x_1\rrbracket\alpha_3\in \alpha_1 H^{s_2+s_3-1,t_2+t_3}(U) +
H^{s_1+s_2-1,t_1+t_2}(U)\,\alpha_3=\Ind\langle\alpha_1,\alpha_2,\alpha_3\rangle$. By
$\llbracket x\rrbracket$ we mean the class of $x$ in $H^*(U)$. The group
$\Ind\langle\alpha_1,\alpha_2,\alpha_3\rangle$ is called the \emph{indeterminacy} of
$\langle\alpha_1,\alpha_2,\alpha_3\rangle$. If the indeterminacy is trivial, then
$\langle\alpha_1,\alpha_2,\alpha_3\rangle$ is a single class in $H^**(U)$.

\begin{definition} With notation as above, $\langle\alpha_1,\alpha_2,\alpha_3\rangle\subset
H^{s-1,t}(U)$ is the coset of $\Ind\langle\alpha_1,\alpha_2,\alpha_3\rangle$ represented
by $\bar{a}_1 u_2+\ol{u}_1a_3$.
\end{definition}

Similarly, one can define $k$-fold Massey products $\langle \alpha_1, \dots, \alpha_k
\rangle$ for $k>3$. Let $\alpha_i\in H^{s_i,\,t_i}$ be represented by $a_i$, and let
$s_{i,j}=\sum_{r=i+1}^{j}(s_r -1)$, \,$t_{i,j}=\sum_{r=i+1}^{j} t_r$,\, $s=s_{1,k}+2$,
and $t= t_{1,k}$.

\begin{definition}
A collection of cochains $\mathcal{A}=\bigl\{a_{i,j}\in
U^{s_{i,j}+1,\,t_{i,j}}\bigr\}_{\scriptscriptstyle{0\leq i < j \leq k} }$ is said to be a
\emph{defining system} for the $k$-fold Massey product $\langle \alpha_1, \dots, \alpha_k
\rangle$ if
\begin{align*}
a_{i-1,i} & =a_i \text{ for } i=1, \dots, k,\\
\delta a_{i,j}& =\sum_{i<r<j}\overline{a}_{i,r} a_{r,j}.\\
\intertext{The cocycle $c\bigl(\mathcal{A}\bigr) \in U^{s,t}$ defined by} c(\mathcal{A})&
= \sum_{r=1}^{k-1}\overline{a}_{1,r} a_{r,k},
\end{align*}
is called the \emph{related cocycle} of $\mathcal{A}$.  The $k$-fold Massey product
$\langle \alpha_1, \dots, \alpha_k \rangle\subset H^{s,t}(U)$ consists of all classes
$\left\llbracket c\bigl(\mathcal{A}\bigr)\right\rrbracket$, where $\mathcal{A}$ ranges
over all defining systems of $\langle \alpha_1, \dots, \alpha_k \rangle$.
\end{definition}


It can be shown that the set $\langle \alpha_1, \dots, \alpha_k \rangle$ does not depend
on the choice of representatives of the $\alpha_i$. The $k$-fold product is then a
cohomology operation of $k$ variables that is defined if all the lower products $\langle
\alpha_i,\dots, \alpha_j\rangle$ for $0<j-i<k-1$ are defined and contain zero.  Here the
double product $\langle \alpha_i,\alpha_{i+1}\rangle$ is understood to be the regular
product $\alpha_i\alpha_{i+1}$. The indeterminacy  can get quite complicated, however
\cite[A1.4.4]{Ravenel}. Kraines has constructed a `restricted' $k$-fold Massey product on
the cohomology of topological spaces which has no indeterminacy \cite[Section
3]{Kraines}. This construction can be carried over to the Hopf algebra setting
\cite[Section 3]{Kochman}.
If $\alpha_i=\alpha$ for all $i$ and $\alpha^2=0$, one can define the \emph{symmetric}
Massey product $\langle\alpha\rangle^k$.

\begin{definition} Let $\alpha\in H^{s,t} (U)$ be a class such that $\alpha^2=0$.
A collection of cocycles
 $\mathcal{A}=\bigl\{a_i\in U^{i(s-1)+1,it}\bigr\}_{i=1}^{k-1}$  is called a \emph{defining system} for the
symmetric  Massey product $\langle\alpha\rangle^k$ if
\begin{align*}
a_1 & \text{ represents } \alpha,\\
\delta a_i& = \sum_{r=1}^{i-1}\overline{a}_r a_{i-r} \text{ for } 2\leq i\leq r-1.\\
\intertext{The cocycle $c\bigl(\mathcal{A}\bigr)\in U^{k(s-1)+2, kt}$ defined by}
c\bigl(\mathcal{A}\bigr)& =\sum_{r=1}^{k-1}\overline{a}_r a_{k-r} \in U^{k(s-1)+2, kt}
\end{align*}
is called the \emph{related cocycle} of $\mathcal{A}$.  The symmetric product
$\langle\alpha\rangle^k$ is then the set of all classes $u\in H^{k(s-1)+2,kt}(U)$ which
are represented by some $ c\bigl(\mathcal{A}\bigr)$, where $\mathcal{A}$ is a defining
system for $\langle\alpha\rangle^k$. If $\langle\alpha\rangle^k$ is defined and equal to
a single homology class then $\langle\alpha\rangle^k$ is said to be defined with
\emph{zero indeterminacy}.
\end{definition}

\begin{lemma} Let $p$ be odd, and let $A$ be a Hopf algebra over $\f$.
\be
\item If $\alpha\in H^{2s+1,t}(A)$, then $\langle\alpha\rangle^p$ is defined with
zero indeterminacy \cite[Lemma 14]{Kochman}.
\item If $\alpha\in H^{1,t}(A)$, then $\langle\alpha\rangle^p=-\bp(x)$ \cite[Remark 11.11]{May}.
\ee
\end{lemma}

\subsection{Calculations}\label{calculations}

Let
\begin{equation}\label{truncated polynomials}
\begin{aligned}
A &=\f[x]\big/\bigl(x^{p^n}\bigr), \quad &\deg x &=m \text{ even, } n\geq 1,\quad p
\text{ odd,\
}\\
A &=\mathbb{F}_2[x]\big/(x^{2^n}),\quad &\deg x &= m,\, n>1,\quad p=2.
\end{aligned}
\end{equation}
The augmentation in both cases is defined by $\epsilon (x)=0$. It is known that the
cohomology of such a truncated polynomial algebra is the tensor product of an exterior
algebra on a one-dimensional class and a polynomial algebra on a two dimensional class.
We would like to identify specific generators in the cobar complex $C^*(A)=\Tens\bl
I(A)^*\br$.

Choose for $I(A)$ the natural basis $\{x^i, \ 0 < i < p^n \}$. Let $\{y_i \}$ be the dual
basis for $I^*$. Then
\begin{equation}\label{a}
\delta \bigl[y_k\bigr]= \sum_{\substack{i+j=k\\ 0<i<k}} \bigl[y_i \vert y_{j} \bigr].
\end{equation}

Clearly, in $C^{1,2m}(A)$,\, $z = [y_1]$ is a cocycle. Let $e= \sum_{i=1}^{p^n -1} [y_i
\vert y_{p^n-i}]$. A routine calculation shows that $e$ is a nonbounding cocycle. By
abuse of notation, we can thus take $z$ and $e$ to be the generators of $H^*(A)$.

\begin{lemma}\label{Massey product for u} Let $A$ be a truncated polynomial
algebra as in \eqref{truncated polynomials}. Then $H^{*}(A) \simeq
\Lambda(z)\otimes\f[e]$, with $z\in H^{1,m}$, $e\in H^{2,mp^n}$ as defined above.
Moreover, the symmetric Massey product ${\langle z\rangle}^{p^n}=e$ with zero
indeterminacy.
\end{lemma}

\begin{proof}
Choosing $a_i=[y_i]$ gives a defining system for $\langle z\rangle^{p^n}$. Its related
cocycle is $\sum_{r=1}^{p^n-1}\bigl[y_r\vert y_{p^n-r}\bigr]$, which represents $e$.
There is no indeterminacy because for each $i$, there is no other choice of $a_i$.
\end{proof}
\begin{remark}\label{remark} If $p=2$, $n=1$, then $A=\mathbb{F}_2[x]/(x^2)\simeq \Lambda[x]$, and
$H^*(A)\simeq \mathbb{F}_2[z]$ with $z=[x^*]\in H^{1,\,\deg x}$. 
\end{remark}

Summarizing,

\begin{theorem}\label{non-zero Steenrod operations} Let $A=\f[x]/(x^{p^n})$,
then the following Steenrod operations and symmetric Massey products are present in
$H^*(A)$:
\begin{enumerate}
\item $\langle z\rangle^{p^n}=e$,\, $\pti{1} e=e^p$ for $p$ odd,
\item $\bp z=-e$ for $p$ odd, $n=1$,
\item $\langle z\rangle^{2^n}=e$,\, $\sqit{2}(e)=e^2$ for $p=2$, $n\geq 2$,
\item $\sqit{1}(z)=z^2$ for $p=2$, $n=1$.
\end{enumerate}
\end{theorem}
\begin{proof} These are consequences of the axioms in Theorem \ref{steenrod} and of the
results just stated.
\end{proof}

\bibliographystyle{amsplain}
\providecommand{\bysame}{\leavevmode\hbox to3em{\hrulefill}\thinspace}
\providecommand{\MR}{\relax\ifhmode\unskip\space\fi MR }
\providecommand{\MRhref}[2]{%
  \href{http://www.ams.org/mathscinet-getitem?mr=#1}{#2}
}
\providecommand{\href}[2]{#2}

\end{document}